\documentstyle[11pt]{article}

\begin{document}
\title{Hilbert $C^*$-modules over monotone complete\\ $C^*$-algebras}
\author{{\sc Michael Frank}}
\maketitle

\renewcommand{\baselinestretch}{1.2}
\begin{abstract}
   The aim of the present paper is to describe self-duality and $C^*$-
   reflexi\-vity of Hilbert {\bf A}-mo\-du\-les $\cal M$ over monotone
   complete $C^*$-al\-ge\-bras {\bf A} by the completeness of the unit ball
   of $\cal M$ with respect to two types of convergence being defined, and by
   a structural criterion. The derived results generalize earlier results of
   {\sc H.~Widom} [Duke Math.~J.~23, 309-324, MR 17 \# 1228] and
   {\sc W.~L.~Paschke} [Trans. Amer.~Math.~Soc.~182, 443-468, MR 50
   \# 8087, Canadian J.~Math.~26, 1272-1280, MR 57 \# 10433]. For Hilbert
   $C^*$-modules over commutative $AW^*$-algebras the equivalence
   of the self-duality property and of the Kaplansky-Hilbert property is
   reproved, (cf. {\sc M.~Ozawa} [J.~Math.~Soc.~Japan 36, 589-609, MR
   85m:46068] ). Especially, one derives that for a $C^*$-algebra {\bf A}
   the {\bf A}-valued inner pro\-duct of  \linebreak[4]
   every Hilbert {\bf A}-module $\cal M$
   can be continued to an {\bf A}-valued inner pro\-duct on it's {\bf A}-dual
   Banach {\bf A}-module $\cal M{\rm '}$ turning $\cal M{\rm '}$ to a
   self-dual Hilbert {\bf A}-module if and only if {\bf A} is monotone
   complete (or, equivalently, additively complete) generalizing a result of
   {\sc M.~Hamana} [Internat.~J.~Math.~3(1992), 185-204].
   A classification of countably generated self-dual Hilbert
   {\bf A}-modules over monotone complete $C^*$-algebras {\bf A} is
   established. The set of all bounded module operators ${\bf End}_A (\cal M)$
   on self-dual Hilbert {\bf A}-modules $\cal M$ over monotone complete
   $C^*$-algebras {\bf A} is proved again to be a monotone complete
   $C^*$-algebra. Applying these results a {\sc Weyl}-{\sc Berg} type theorem
   is proved.
\end{abstract}

   \newpage

   Recently, some progress was made in the theory of operator valued
   weights and conditional expectations of finite index on $W^*$-algebras
   \cite{Bai} , and in a generalized {\sc Tomita}-{\sc Takesaki} theory for
   embeddable $AW^*$-algebras \cite{Fr2} using self-dual Hilbert $W^*$-modules
   and Kaplansky-Hilbert modules as the essential technical tool. However, as
   almost everywhere in the literature the theory was developed only as far as
   necessary for the solution of the main problem of the authors.
   Since the tool "Hilbert $C^*$-module" has been used by more and more
   mathematicians, previously in operator theory, it might be quite useful to
   sharpen the tool for getting deeper results in it's areas of application.
   All the more there are only a few such well-known basic publications on
   this subject like \cite{Kp2}, \cite{Pa1,Pa2}, \cite{Rie}, \cite{Ks},
   \cite{Lin}, \cite{Ham92}, \cite{La}, for example. The present paper is an
   attempt to make a next step.

   The concrete cause for the investigations below are two problems formulated
   by \linebreak[4] {\sc W.~L.~Paschke} (cf.~\cite{Pa1,Pa2}) , and the problem
   of finding a purely $C^*$-algebraical approach to the diagonalization of
   operators up to a suitable small rest (cf.~\cite{Zh1,Zh2}).

   {\sc W.~L.~Paschke} raised the question for which $C^*$-algebras {\bf A} the
   {\bf A}-valued inner product on an arbitrary Hilbert {\bf A}-module
   $\cal M$ can be lifted to an {\bf A}-valued inner product on its
   {\bf A}-dual Banach {\bf A}-module $\cal M{\rm '}$ turning $\cal M{\rm '}$
   into a self-dual Hilbert {\bf A}-module. By his construction one can find
   such a lifting in the case of {\bf A} being a $W^*$-algebra or a
   commutative $AW^*$-algebra. On the other hand he showed that {\bf A} has
   to be at least an $AW^*$-algebra, cf.~\cite[ Th.3.2.]{Pa1},
   \cite[ Prop.1]{Pa2}. Refering to {\sc M.~Hamana} \cite[Th. 2.2]{Ham92} and
   {\sc H.~Lin} \cite[Lemma 3.7]{Lin} we resolve this problem for monotone
   complete $C^*$-algebras affirmatively giving a general construction based
   on order convergence. Moreover, glueing together results of {\sc M.~Hamana}
   \cite{Ham92} and {\sc K.~Sait{\^o}}, {\sc J.~D.~M.~Wright} \cite{Sai/Wri}
   we can show that this is the general solution of {\sc W.~L.~Paschke}'s
   problem. However, the central problem of the $AW^*$-theory, whether all
   $AW^*$-algebras are monotone complete or not, is still open despite
   recent encouraging results (\cite{Chr/Ped,Sai/Wri}).

   The second problem under consideration is to describe the inner
   structure of self-dual Hilbert {\bf A}-modules $\cal M$ over $C^*$-algebras
   {\bf A}. {\sc W.~L.~Paschke} has found very nice criteria on self-duality and
   $C^*$-reflexivity in the case of {\bf A} being a $W^*$-algebra or
   a commutative $AW^*$-algebra. These ideas were extended by other
   authors. (Cf.~\cite{Pa1},\cite{Fr1},\cite{Bai}). Moreover, the author was
   able to show that the self-duality of an arbitrary Hilbert {\bf A}-module
   $\cal M$ over a $C^*$-algebra {\bf A} does not depend on the
   {\it properties} of the concrete given {\bf A}-valued inner product on
   $\cal M$ realizing the self-duality, but only depends on the {\it existence}
   of such an {\bf A}-valued inner product on the Banach {\bf A}-module
   $\cal M$ inducing an equivalent norm to the given one, \cite{Fr1}.
   Therefore, one has the possibility to characterize self-duality of
   arbitrary Hilbert $C^*$-modules by their inner structure. We
   give such a characterization in the case of {\bf A} being a monotone
   complete $C^*$-algebra. On the way the interrelation is shown between the
   theory of Kaplansky-Hilbert modules and the theory of self-dual Hilbert
   {\bf A}-modules over commutative $AW^*$-algebras,
   cf.~\cite[ Th.5.4]{Oz1},\cite{ESW},\cite{Kp2},\cite{Oz3},\cite{Wd}.

   These results allow us to prove some generalized {\sc Weyl}-{\sc Berg}
   type theorems for boun\-ded, A-linear, normal operators on self-dual Hilbert
   {\bf A}-modules with countably generated {\bf A}-pre-dual Hilbert
   {\bf A}-modules, where {\bf A} is assumed to be monotone complete and to
   have a special approximation property  (*). An example shows that a
   {\sc Weyl} type decomposition can fail if the $C^*$-algebra {\bf A} does
   not have property (*).

   Let us remark that the theory of operator valued weights and conditional
   expectations of finite index on $W^*$-algebras in the approach of
   {\sc M.~Baillet}, {\sc Y.~Denizeau} and \linebreak[4] {\sc J.-F.~Havet}
   \cite{Bai} can be ge\-ne\-ra\-lized straightforward to the case of
   monotone complete $C^*$-algebras using the results of the present paper.
   But this appears elsewhere, cf.~\cite{Fr4}. The main content of the
   present paper was previously circulated as a preprint, \cite{Fr3}.

   \medskip
   The present paper is organized as follows:

   \smallskip
   After a section about basic definitions and facts the second one is
   concerned with two types of convergence in Hilbert {\bf A}-modules over
   monotone complete $C^*$-algebras {\bf A} and their properties. In section
   three a spectral decomposition theorem for normal elements of normal
   $AW^*$-algebras and a polar decomposition theorem for arbitrary
   $AW^*$-algebras are proved. Two criteria on self-duality and
   $C^*$-reflexivity of Hilbert $C^*$-modules are established in section four
   and five. They have several consequences of interest: For every Hilbert
   {\bf A}-module $\{ \cal M, \langle .,. \rangle \}$ over a monotone complete
   $C^*$-algebra {\bf A} the {\bf A}-valued inner product $\langle .,. \rangle$
   can be continued to an {\bf A}-valued inner product on its {\bf A}-dual
   Banach {\bf A}-module $\cal M{\rm '}$ turning $\cal M{\rm '}$ into a
   self-dual Hilbert {\bf A}-module, (section four). One obtains a
   classification of self-dual, countably generated Hilbert {\bf A}-modules
   over monotone complete $C^*$-algebras {\bf A}, (section six). The
   $C^*$-algebra ${\bf End}_A(\cal M)$ of all bounded, {\bf A}-linear
   operators on self-dual Hilbert {\bf A}-modules $\cal M$ over monotone
   complete $C^*$-algebras {\bf A} turns out to be monotone complete, again.
   Hence, spectral and polar decomposition are possible inside
   ${\bf End}_A(\cal M)$, and one can consider some kind of {\sc von Neumann}
   representations of monotone complete $C^*$-algebras on such Hilbert
   $C^*$-modules (section seven). This last result gives us the chance to
   prove a generalized {\sc Weyl}-{\sc Berg} type theorem for the pair
   $\{ {\bf End}_A(\cal M),{\bf K}_A(\cal M)\}$ with {\bf A} being a monotone
   complete $C^*$-algebra with the approximation property (*). (Here
   ${\bf K}_A(\cal M)$ denotes the set of all bounded, "compact" operators on
   $\cal M$.) This is explained in section seven.

\section{Preliminaries and basic facts}
   Before we start considerations let us recall some definitions and fix the
   notations. We consider {\it $AW^*$-algebras}, i.e. $C^*$-algebras for which
   the following two conditions are satisfied, (\cite{Kp1}):

   \smallskip \noindent
   (a) In the partially ordered set of projections every set of pairwise
   orthogonal projections has a least upper bound.

   \noindent
   (b) Every maximal commutative self-adjoint subalgebra is generated by
   its projections, i.e., it is equal to the smallest closed *-subalgebra
   containing its projections.

   \smallskip
   A $C^*$-algebra is said to be {\it monotone complete} if and only if every
   bounded increa\-singly directed net $ \{ a_\alpha :\alpha \in I \} $ of
   self-adjoint elements of it has a least upper bound \linebreak[4]
   $ a= \sup \{ a_\alpha: \alpha \in I \} $ in it. Every monotone complete
   $C^*$-algebra is an $AW^*$-algebra. At present it is unknown whether every
   $AW^*$-algebra is monotone complete or not, cf.~\cite{ESW}. But commutative
   $AW^*$-algebras are monotone complete, (\cite{Az}). More examples of
   monotone complete $AW^*$-algebras can be found at \cite[Th.7, Th.8]{Kp2},
   \cite[Th.1]{Pa1}, \cite{ESW}. An $AW^*$-algebra is called {\it
   normal} if every increasingly directed net of projections $\{p_{\alpha}:
   \alpha \in I \}$ with supremum $p$ inside the net of all projections
   possesses a supremum with respect to the net of all self-adjoint elements
   (being equal to $p$ automatically in the case of it's existence),
   cf.~\cite{Wright:80}. A $C^*$-algebra is {\it additively complete} if every
   norm-bounded sum of positive elements has a least upper bound.

   If {\bf A} is a monotone complete $C^*$-algebra then one can define order
   convergence on {\bf A}.
   For the commutative case this was done by {\sc H.~Widom} \cite{Wd}
   (cf.~\cite{Azar}), and for the general case by {\sc R.~V.~Kadison} and
   {\sc G.~K.~Pedersen} \cite{KP} who defined the so-called Kadison-Pedersen
   arrow. {\sc M.~Hamana} \cite[p.260]{Hm2} modified the latter notion to get
   a general notion of order convergence.
   A net $ \{ a_\alpha : \alpha \in I \} $ of elements of {\bf A} {\it
   converges to an element $ a \in${\bf A} in order} if and only if there are
   bounded nets $ \{ a^{(k)}_\alpha : \alpha \in I \} $ and $ \{
   b^{(k)}_\alpha : \alpha \in I \} $ of self-adjoint elements of {\bf A} and
   self-adjoint elements $ a^{(k)} \in${\bf A}, $k=1,2,3,4 $, such that

   \hspace{0.5cm} (i)  $ \: 0 \le a^{(k)}_\alpha - a^{(k)} \le b^{(k)}_\alpha
                       ,\: k=1,2,3,4, \alpha \in I $,

   \hspace{0.5cm} (ii)  $ \, \{ b^{(k)}_\alpha : \alpha \in I \} $  is
                        decreasingly directed and has greatest lower bound
                        zero,

   \hspace{0.5cm} (iii) $ \sum_{k=1}^4 ({\bf i})^{k} a^{(k)}_\alpha =
                        a_\alpha $ for every $ \alpha \in I $ ,
                        $ \sum_{k=1}^{4} ({\bf i})^{k} a^{(k)} = a $ ,
                        (where  $ {\bf i}=\sqrt{-1} $).

   \noindent
   We denote this type of convergence by  LIM$\{ a_\alpha : \alpha \in I \}=a$. By
   \cite[p.260]{Hm1} the order limit of $\{ a_\alpha : \alpha \in I \}$  does
   not depend on the special choice of the nets $\{ a^{(k)}_\alpha : \alpha
   \in I \}, \{ b^{(k)}_\alpha : \alpha \in I \}$ and of the elements
   $a^{(k)} ,\: k=1,2,3,4 $.
   If {\bf A} is a commutative $AW^*$-algebra then the order convergence
   defined above is equivalent to the classical order convergence in {\bf A}
   as it was defined by {\sc H.~Widom} \cite{Wd} earlier. Order convergence
   has the following properties, cf.~\cite[Lemma 1.2]{Hm2}:

   \smallskip
   If  LIM$\{ a_\alpha : \alpha \in I \} = a$ , LIM$\{ b_\beta : \beta \in J
   \} = b $ then

   \hspace{0.3cm} (i) $ \,\; $LIM$\{ a_\alpha + b_\beta : \alpha \in I, \beta
                      \in J \} = a+b $,

   \hspace{0.3cm} (ii) $\;$LIM$\{ x a_\alpha y : \alpha \in I \} = xay $ for
                       every $ x,y \in A $,

   \hspace{0.3cm} (iii)  LIM$\{ a_\alpha b_\beta : \alpha \in I, \beta \in J
                         \} = ab $,

   \hspace{0.3cm} (iv)  $\; a_\alpha \le b_\alpha $ for every $ \alpha \in I $
                        implies $ a \le b $,

   \hspace{0.3cm} (v)  $ \,\; \| a \|_A \le \limsup \{ \| a_\alpha \|_A :
                       \alpha \in I \} $.

   \smallskip
   Throughout this paper we denote the $C^*$-norm of {\bf A} by $ \| . \|_A $.
   The self-adjoint part of {\bf A} is denoted by $ {\bf A}_h $ and the
   positive cone of {\bf A} by $ {\bf A}^+_h $. The centre of {\bf A} has the
   denotation ${\bf Z}({\bf A})$.

   \medskip
   Now some facts about Hilbert $C^*$-modules. We make the convention that all
   $C^*$-modules of the present paper are left modules by definition. A
   {\it pre-Hilbert {\bf A}-module over a fixed $C^*$-algebra} {\bf A} is
   an {\bf A}-module $\cal M$ equipped with an {\bf A}-valued, non-degenerate
   mapping \linebreak[4] $ \langle .,. \rangle : \cal M \times \cal M
   \longrightarrow {\bf A} $ being {\bf A}-linear in the first argument and
   conjugate-{\bf A}-linear in the second, and satisfying $\langle x,x
   \rangle \in {\bf A}_h^+$ for every $x \in \cal M$. The map $ \langle .,.
   \rangle $ is called {\it the {\bf A}-valued inner product on} $\cal M$.
   A pre-Hilbert {\bf A}-module $ \{ \cal M, \langle .,. \rangle \} $ is
   {\it Hilbert} if and only if it is complete with respect to the norm
   $ \| . \| = \| \langle .,. \rangle \|^{1/2}_A $. We always suppose that the
   linear structures of  {\bf A} and $\cal M$ are compatible. Denote by
   $ \langle \cal M , \cal M \rangle $ the norm closed linear hull of the
   range of the map $ \langle .,. \rangle $ in {\bf A}. A Hilbert
   {\bf A}-module $ \{ \cal M , \langle .,. \rangle \} $ over a $C^*$-algebra
   {\bf A} is said to be {\it self-dual} if and only if every bounded module
   map $ r: \cal M \longrightarrow {\bf A} $ is of the form $ \langle . ,
   a_r \rangle $ for some element $ a_r \in \cal M $. The set of all
   bounded module maps $ r: \cal M \longrightarrow {\bf A} $ is denoted by
   $\cal M{\rm '}$ . It is a Banach {\bf A}-module. A Hilbert {\bf A}-module
   is said to be {\it $C^*$-reflexive} (or {\it {\bf A}-reflexive}) if and
   only if the map $ \Omega $ being defined by the formula
   \[
   \Omega (x)[r]=r(x) \quad {\rm for}\: {\rm each}\:  x \in \cal M \:,\quad
   {\rm every} \: {\it r} \in \cal M' ,
   \]
   is a surjective module mapping of $\cal M$ onto the Banach {\bf A}-module
   $\cal M{\rm ''}$, where $\cal M{\rm ''}$ consists of all bounded module
   maps from $\cal M{\rm '}$ to {\bf A}. A Hilbert C*-module $\cal M$ is
   {\it countably generated} if there exists a countable set of elements of
   $\cal M$ the set of finite C*-linear combinations of which being norm-dense
   in $\cal M$.
   For more basic facts about Hilbert
   $C^*$-modules we refer to \cite{Pa1,Pa2,La}.

\section{Two types of convergence}

   {\bf Definition 2.1 :} (cf.~\cite[\S2]{Wr2}, \cite[\S1.1]{Wd},
   \cite[Def.3.1]{Fr1}, \cite{Azar})
   Let {\bf A} be a monotone complete $C^*$-algebra, $\{\cal M,\langle .,.
   \rangle\}$ be a pre-Hilbert {\bf A}-module and $I$ be a net. A norm-bounded
   set \linebreak[4] $\{x_\alpha : \alpha \in I\}$ of elements of $\cal M$
   is {\it fundamental in the sense of} $\tau_1^0$-{\it convergence} (in short:
   $\tau_1^0$-{\it fundamental}) if and only if the limits
   \[
   {\rm LIM}\{ \langle x_\alpha - x_\beta , x_\alpha - x_\beta \rangle : \alpha \in I\}
   \]
   exist for every $\beta \in I$, and the limit
   \[
   {\rm LIM}\{ {\rm LIM}\{ \langle x_\alpha - x_\beta , x_\alpha - x_\beta
   \rangle :\alpha \in I \} : \beta \in I \}
   \]
   exists, too, and equals zero. A norm-bounded set $\{ x_\alpha : \alpha \in
   I \}$ of elements of $\cal M$ {\it has the} $\tau_1^0${\it -limit} $x$
   {\it in} $\cal M$ if and only if the limit
   \[
   {\rm LIM}\{ \langle x_\alpha - x , x_\alpha -x \rangle : \alpha \in I \}
   \]
   exists and equals zero. In this case one writes
   \[
   \tau_1^0-\lim \{ x_\alpha : \alpha \in I \} = x,
   \]
   and one says that $\{ x_\alpha : \alpha \in I \}$ $\tau_1^0${\it
   -converges to $x$ inside} $\cal M$.

   \medskip
   {\bf Lemma 2.2 :} {\it Let} {\bf A} {\it be a monotone complete} $C^*${\it
   -algebra and} $\{ \cal M, \langle .,. \rangle\}$ {\it a Hilbert}
   {\bf A}{\it -module. If } $\{ x_\alpha : \alpha \in I \}$ {\it is an
   norm-bounded set of elements of} $\cal M$ {\it  indexed by a net} $I$ {\it
   and possessing a} $\: \tau_1^0${\it -limit} $\: x$ {\it in} $\cal M$, {\it
   then the following statements are true:}

   (i)  $\;$ {\it The} $\tau_1^0${\it -limit} $x \in \cal M$ {\it is unique.}

   (ii) {\it The set} $\{ x_\alpha : \alpha \in I  \}$ {\it is} $\tau_1^0${\it
   -fundamental.}

   (iii) {\it If} $\| x_\alpha \| \le N$ {\it for every} $\alpha \in I$ {\it
   and a real number} $N$ {\it then}  $\| x \| \le N$.

   (iv) {\it For every} $a \in${\bf A} {\it the equality} $\tau_1^0-\lim
   \{ ax_\alpha : \alpha \in I \} = ax$ {\it is satisfied.}

   (v) {\it For every norm-bounded set} $\{ y_\beta : \beta \in J \}$
   {\it (not necessarily distinct from} \linebreak[4] $\{ x_\alpha : \alpha
   \in I \})$ {\it of elements of} $\cal M$ {\it indexed by a net} $J$ {\it
   and possessing a} $\tau_1^0${\it -limit} $y$ {\it in} $\cal M$ {\it the
   equality} $\tau_1^0-\lim \{ x_\alpha + y_\beta : \alpha \in I, \beta \in J
   \} = x+y$ {\it holds.}

   \medskip
   P r o o f:  To prove the first fact suppose the existence of two
   $\tau_1^0$-limits $x_1, x_2 \in \cal M$ of the norm-bounded
   net $\{ x_\alpha : \alpha \in I \} $. The inequality
   \[
   0 \le \langle x_1-x_2 , x_1-x_2 \rangle \le 2(\langle x_\alpha-x_1 ,
   x_\alpha-x_1 \rangle + \langle x_\alpha-x_2 , x_\alpha-x_2 \rangle)
   \]
   holds for every $\alpha \in I$. If one calculates the order limit of the
   right side one obtains the equa\-lity $x_1 = x_2$. Furthermore, the fact
   (ii) can be derived from the inequality
   \[
   0 \le \langle x_\alpha-x_\beta , x_\alpha-x_\beta \rangle \le
   2(\langle x_\alpha-x , x_\alpha-x \rangle + \langle x_\beta-x , x_\beta-x
   \rangle)
   \]
   being satisfied for every $\alpha,\beta \in I$, and from the
   $\tau_1^0$-convergence of the set $\{ x_\alpha : \alpha \in I \}$ to $x
   \in \cal M$. Similarly, (v) follows from the inequality
   \[
   0 \le \langle x+y-x_\alpha-y_\beta , x+y-x_\alpha-y_\beta \rangle \le
   2(\langle x-x_\alpha , x-x_\alpha \rangle + \langle y-y_\beta , y-y_\beta
   \rangle)
   \]
   being satisfied for every $\alpha \in I$, $\beta \in J$, and from the
   $\tau_1^0$-convergence of the sets \linebreak[4] $\{ x_\alpha : \alpha \in
   I \}, \{ y_\beta : \beta \in J \} \: {\rm to} \: x,y \in \cal M$,
   respectively. To show (iii) notice that by (v) there exists the order limit
   of the net $\{ \langle x_\alpha , x_\alpha \rangle : \alpha \in I \} $ in
   {\bf A} being equal to $\langle x,x \rangle$. Therefore,
   \[
   \|x\| = \| \langle  x,x \rangle \|_A^{1/2} \le
   \limsup \{ \| \langle x_\alpha , x_\alpha \rangle \|_A^{1/2} : \alpha \in
   I \} \le N.
   \]
   Statement (iv) is derived from the fact that the inequality $0 \le
   \langle x,x \rangle \le \langle y,y \rangle$ implies the inequality $ 0 \le
   a \langle x,x \rangle a^* \le a\langle y,y \rangle a^*$ for every $x,y \in
   \cal M$, every $a \in${\bf A}, (\cite[Prop.2.2.13]{Bra}).    $\: \bullet$

   \medskip
   {\bf Definition 2.3 :} (cf.~\cite[\S 1.1]{Wd}, \cite[Def.3.1]{Fr1})
   Let {\bf A} be a monotone complete $C^*$-algebra and let $\{ \cal M,
   \langle .,. \rangle \}$ be a pre-Hilbert {\bf A}-module and $I$ be a net.
   A norm-bounded set $\{ x_\alpha : \alpha \in I \}$ of elements of $\cal M$
   is {\it fundamental in the sense of} $\tau_2^0$-{\it convergence} (in short:
   $\tau_2^0$-{\it fundamental}) if and only if the limits
   \[
   {\rm LIM}\{ \langle y,x_\alpha-x_\beta \rangle : \alpha \in I \}
   \]
   exist for every $y \in \cal M$, for every $\beta \in I$, and the limit
   \[
   {\rm LIM}\{ {\rm LIM}\{ \langle y,x_\alpha-x_\beta \rangle : \alpha \in I
   \} : \beta \in I \}
   \]
   exists, too, and equals zero. A norm-bounded set $\{ x_\alpha : \alpha \in
   I\}$ of elements of $\cal M$ {\it has the} $\tau_2^0${\it -limit} $\: x$
   {\it in} $\cal M$ if and only if for every $y \in \cal M$ the limit
   \[
   {\rm LIM}\{ \langle y,x_\alpha-x \rangle : \alpha \in I \}
   \]
   exists and equals zero. In this case one writes
   \[
   \tau_2^0-\lim \{ x_\alpha : \alpha \in I \} = x,
   \]
   and one says that $\{ x_\alpha : \alpha \in I \}$ $\tau_2^0${\it -converges
   to $x$ inside} $\cal M$.

   \medskip
   {\bf Lemma 2.4 :} {\it Let} {\bf A} {\it be a monotone complete} $C^*${\it
   -algebra and} $\{ \cal M, \langle .,. \rangle \}$ {\it be a Hilbert}
   {\bf A}{\it -module. If} $\{ x_\alpha : \alpha \in I \}$ {\it is a
   norm-bounded set of elements of} $\cal M$ {\it indexed by a net I and
   possessing a} $\tau_2^0${\it -limit} $\: x$ {\it in} $\cal M$, {\it then
   the following statements are true:}

   (i) $\;$ {\it The} $\tau_2^0${\it -limit} $x \in \cal M$ {\it is unique.}

   (ii) {\it The set} $\{ x_\alpha : \alpha \in I \}$ {\it is} $\tau_2^0${\it
   -fundamental.}

   (iii) {\it If} $\| x_\alpha \| \le N$ {\it for every} $\alpha \in I$ {\it
   and a real number} $\: N$ {\it then} $\| x \| \le N$.

   (iv) {\it For every} $a \in${\bf A} {\it the equality} $\tau_2^0-\lim \{
   ax_\alpha : \alpha \in I \} = ax$ {\it is satisfied.}

   (v) {\it For every norm-bounded set} $\{ y_\beta : \beta \in J \}$ {\it
   (not necessarily distinct from} \linebreak[4] $\{ x_\alpha : \alpha \in I
   \})$ {\it of elements of} $\cal M$ {\it indexed by a net} $\: J$ {\it and
   possessing a} $\tau_2^0${\it -limit} $\:y$ {\it in} $\cal M$ {\it the
   equality} $\tau_2^0-\lim \{ x_\alpha + y_\beta : \alpha \in I, \beta  \in
   J \} = x+y$ {\it holds.}

   \medskip
   P r o o f: The facts (i), (ii) and (v) are obvious.
   To show (iv) one has only to consider the equality

   \begin{eqnarray*}
      {\rm LIM} \{ \langle y,ax_\alpha \rangle : \alpha \in I \} & = &
                  {\rm LIM} \{ \langle y,x_\alpha \rangle a^* : \alpha \in I
                  \} \\
      & = &  ( {\rm LIM} \{ \langle y,x_\alpha \rangle : \alpha \in I \})
             a^* \\
      & = & \langle y,x \rangle a^* \\
      & = & \langle y,ax \rangle
   \end{eqnarray*}

   being valid for every $a \in${\bf A} and every $y \in \cal M$. For the
   proof of (iii) consider the in\-equa\-lity
   \[
   0 \le \| \langle x_\alpha,y \rangle + \langle y,x_\alpha \rangle \| \le
   2 \|y\|\: \|x\| \le 2 \|y\| N
   \]
   being valid for every $y \in \cal M$, every $\alpha \in I$ and some
   real positive  number $\: N$ by assumption. The element
   $(\langle x_\alpha,y \rangle + \langle y,x_\alpha \rangle)$ is self-adjoint
   for every $y \in \cal M$, every $\alpha \in I$. Therefore,
   \[
   -2 N \|y\| 1_A \le \langle x_\alpha,y \rangle + \langle y,x_\alpha \rangle
   \le 2 N \|y\| 1_A
   \]
   for every $\alpha \in I$, every $y \in \cal M$. Replacing $\: y$ by $\: x$
   and taking the order limit of the central expression statement (iii)
   turns out. $\: \bullet$

   \medskip
   {\bf Lemma 2.5 :} {\it Let} {\bf A} {\it be a monotone complete} $C^*${\it
   -algebra,} $\{ \cal M,\langle .,. \rangle \}$ {\it be a Hilbert}
   {\bf A}{\it -module and} $I$ {\it be a net. If a norm-bounded set} $\{
   x_\alpha : \alpha \in I \}$ {\it of elements of} $\cal M$ {\it possesses a}
   $\tau_1^0${\it -limit} $x$ {\it in} $\cal M$ {\it then that element} $x$
   {\it is the } $\tau_2^0${\it -limit of the net} $\{ x_\alpha : \alpha \in
   I \}$ {\it in} $\cal M$, {\it too}.

   \medskip
   P r o o f:  The statement above follows from the following equation
   based on the pola\-ri\-zation formula:

   \begin{eqnarray*}
        {\rm LIM} \{ \langle x-x_\alpha,y \rangle : \alpha \in I \} & = &
        \frac{1}{4} \; {\rm LIM} \left\{ \sum_{k=0}^3 {\bf i}^k \langle
        x-x_\alpha + {\bf i}^k y,x-x_\alpha + {\bf i}^k y \rangle \right\} \\
            & = &  \frac{1}{4} \; \sum_{k=0}^3  {\bf i}^k {\rm LIM} \{ \langle
            x-x_\alpha + {\bf i}^k y,x-x_\alpha + {\bf i}^k y \rangle \} \\
            & = &  0
   \end{eqnarray*}

   \noindent
   being valid for ${\bf i} = \sqrt{-1}$, every $y \in \cal M$, every
   $\alpha \in I$. Therefore, $x = \tau_2^0-\lim \{ x_\alpha : \alpha \in I
   \}$.      $\: \bullet$

   \medskip
   For example, if one defines an {\bf A}-valued inner product $\langle .,.
   \rangle_A$ by the formula $\langle a,b \rangle :=a b^* , (a,b \in {\bf A
   })$, on a monotone complete  $C^*$-algebra {\bf A} then the
   $\tau_2^0$-convergence induces the order convergence on {\bf A}. On the
   other hand, $\tau_1^0$-convergence induces the order convergence on {\bf A}
   if and only if {\bf A} is commutative. As another example one can consider
   an arbitrary Hilbert space $\cal H$, (where {\bf A} is the set of complex
   numbers). Then $\tau_1^0$-convergence is induced by the norm topology on
   $\cal H$, whereas the $\tau_2^0$-convergence is induced by the weak
   topology on $\cal H$, i.e. the two types of $\tau^0$-convergence do not
   coincide, in general. Moreover, on Hilbert {\bf A}-modules over
   $W^*$-algebras {\bf A} $\tau_1^0$-convergence is induced by a topology
   which is generated by the semi-norms $\{ f(\langle .,. \rangle)^{1/2} :
   f \in {\bf A}_* \}$, and $\tau_2^0$-convergence is induced by a topology
   which is generated by the linear functionals $\{ f(\langle .,x \rangle) :
   f \in {\bf A}_*, x \in \cal M \}$, cf.~\cite{Fr1}.

   \medskip
   {\bf Remark:} Does there exist a topology on monotone complete
   $C^*$-algebras inducing order convergence on them? The answer is negative
   even in the commutative case. A simple example has been constructed by
   {\sc E.~E.~Floyd} \cite{Floyd} in 1955.

   \medskip
   {\bf Proposition 2.6 :} {\it Let} {\bf A} {\it be a monotone complete}
   $C^*${\it -algebra and } $\{ \cal M, \langle .,. \rangle \}$ {\it a
   Hilbert} {\bf A}{\it -module. Consider the standard embedding of} $\cal M$
   {\it into its} {\bf A}-{\it dual Banach} {\bf A}-{\it module}
   ${\cal M}{\rm '}$.
   {\it The linear hull of the} $\tau_1^0${\it -completed unit ball of}
   ${\cal M} \hookrightarrow {\cal M}{\rm '}$ {\it can be identified with a
   Hilbert} {\bf A}{\it -submodule} $\{ {\cal M}^*, \langle .,. \rangle_D \}$
   {\it of the Banach {\bf A}-module ${\cal M}{\rm '}$. The embedding}
   ${\cal M} \hookrightarrow {\cal M}'$ {\it is realized via the mapping}
   $x \in {\cal M} \longrightarrow \langle . , x \rangle_D \in \cal M^*$
   {\it since} $\langle x,y \rangle_D = \langle  x,y \rangle$
   {\it for every} $x,y \in \cal M$.

   \medskip
   P r o o f: By Lemma 2.5 every $\tau_1^0$-fundamental set $\{ x_\alpha :
   \alpha \in I \}$ of $\cal M$ is $\tau_2^0$-fundamental, too. Embedding
   $\cal M$ into ${\cal M}{\rm '}$ in the canonical way and defining
   \[
   r(y) = {\rm LIM} \{ \langle y,x_\alpha \rangle : \alpha \in I \}
   \]
   for arbitrary elements $y \in \cal M$ one finds $r=\tau_2^0-{\rm lim} \{
   x_\alpha : \alpha \in I \}$. Let us denote the union of ${\cal M}
   \hookrightarrow {\cal M}{\rm '}$ with the set of all $\tau_2^0$-limits of
   it's $\tau_1^0$-fundamental subsets by ${\cal M}^*$.

   The next step is to define an {\bf A}-valued inner product on ${\cal M}^*
   \subseteq {\cal M}{\rm '}$. Set:
   \[
   \langle x,y \rangle_D = \langle x,y \rangle \, , \,
   x,y \in {\cal M} \hookrightarrow {\cal M}^* ,
   \]
   \[
   \langle x,r \rangle_D = r(x) \, , \, x \in {\cal M}
   \hookrightarrow {\cal M}^* \, , \, r \in {\cal M}^* \setminus {\cal M} ,
   \]
   This setting is correct because of the definition of $r \in {\cal M}^*
   \setminus {\cal M}$ above. Moreover, since the inequality
   $0 \leq r(y_\gamma -y_\delta)^*r(y_\gamma -y_\delta)
   \leq \| r \|^2 \langle y_\gamma -y_\delta,y_\gamma -y_\delta \rangle$
   is valid for every $\tau_1^0$-fundamental set $\{ y_\gamma : \gamma \in J
   \}$ one defines
   \[
   \langle s,r \rangle_D = {\rm LIM} \{ r(y_\gamma) : \gamma \in J \} \, ,
   \, r,s \in {\cal M}^* \setminus {\cal M} \, , \, s = {\rm LIM} \{ y_\gamma
   : \gamma \in J \} .
   \]
   Note, that the sequence of the two order limits is irrelevant to the result.
   With such a setting define arbitrary values of the {\bf A}-valued inner
   product with respect to the axioms of this structure. By Lemma 2.2, 2.4 and
   2.5 the {\bf A }-valued inner product $\langle .,. \rangle_D$ is
   well-defined.

   Now, we have to show that the $\tau_2^0$-limit $r$ of the
   $\tau_1^0$-fundamental set $\{ x_\alpha : \alpha \in I \}$ of $\cal M$ is
   its $\tau_1^0$-limit, too. Obviously,
   \begin{eqnarray*}
   {\rm LIM} \{ (\langle x_\alpha -r,x_\alpha -r \rangle ) : \alpha \in I \}
   & = &
      {\rm LIM} \{ (\langle x_\alpha ,x_\alpha - x_\beta \rangle -
      r(x_\alpha - x_\beta )^*) : \alpha , \beta \in I \} \\
   & = & {\rm LIM} \{ (\langle x_\alpha ,x_\alpha - x_\beta \rangle -
   \langle x_\beta ,x_\alpha - x_\beta \rangle ) : \alpha , \beta \in I \} \\
   & = & 0 \, .
   \end{eqnarray*}
   Hence, every $\tau_1^0$-fundamental set of $\cal M$ has a $\tau_1^0$-limit
   inside ${\cal M}^*$. The set $\cal M^*$ is obviously
   complete with respect to the norm $\| \langle .,. \rangle_D \|_A^{1/2}$
   by its definition.   $\: \bullet$

\section{ Spectral and polar decomposition inside $AW^*$-algebras}

   We want to show a spectral theorem for normal elements of normal
   $AW^*$-algebras. It gives important informations about normal and,
   especially, self-adjoint elements being necessary for further
   considerations. To formulate the theorem the following definition is
   useful:

   \medskip
   {\bf Definition 3.1 :} (cf. \cite[p.264]{Wr2})
   A measure $m$  on a compact Hausdorff space $X$ with values in the
   self-adjoint part of a monotone complete $C^*$-algebra is called {\it
   quasi-regular} if and only if
   \[
   m(K) = \inf \{ m(U) : U {\rm - open}\:{\rm sets}\:{\rm in} \: X, \: K
   \subseteq U \}
   \]
   for every closed set $K \subseteq X$ . We remark that the condition:
   \[
   m(U) = \sup \{ m(K) : K{\rm  - closed}\:{\rm sets,}\: K \subseteq U \}\:
   {\rm for}\:{\rm every}\:{\rm open}\:{\rm set}\: U \subseteq X
   \]
   is equivalent to quasi-regularity. Further, the measure $m$ is called 
   {\it regular} if and only if
   \[
   m(E) = \inf \{ m(U) : U {\rm - open}\:{\rm sets}\:{\rm in}\: X, E
   \subseteq U \}
   \]
   for every Borel set $E \subseteq X$.

   \medskip
   {\bf Theorem 3.2 :} (cf. \cite[Th.3.1, Th.3.2]{Wr2}
   {\it Let} {\bf A} {\it be a normal} $AW^*${\it -algebra and let} $a \in
   {\bf A}$ {\it be a normal element. Let} ${\bf B} \subseteq {\bf A}$ {\it
   be the commutative} $C^*${\it -subalgebra in} {\bf A} {\it generated by the
   elements} $\{ 1_A, a, a^* \}${\it , and denote by} $\hat{{\bf B}}$ {\it the
   smallest commutative} $AW^*${\it -algebra inside} {\bf A} {\it containing}
   {\bf B} {\it and being monotone complete inside every maximal commutative}
   $C^*${\it -subalgebra {\bf D} of {\bf A} with the property} ${\bf B}
   \subseteq {\bf D}$. {\it Then there exists a unique quasi-regular}
   $\hat{{\bf B}}${\it -valued measure} $m$ {\it on the spectrum} $\sigma(a)
   \subset {\bf C}$ {\it of} $a \in {\bf A}$, {\it the values of which are
   projections in} $\hat{{\bf B}}$ {\it and for which the integral}
   \[
   \int_{\sigma(a)} \lambda \: dm_\lambda \quad = \quad a
   \]
   {\it exists in the sense of order convergence in} $\hat{{\bf B}} \subseteq
   A$.

   \medskip
   P r o o f: By {\sc Gelfand}-{\sc Naimark} the commutative $C^*$-subalgebra
   ${\bf B} \subseteq {\bf A}$ being ge\-ne\-rated by the elements $\{ 1_A, a,
   a^* \}$ is *-isomorphic to the commutative $C^*$-algebra $C( \sigma(a) )$
   of all complex valued continuous functions on the spectrum $\sigma(a)
   \subset {\bf C}$ of $a \in {\bf A}$. Denote this *-isomorphism by $\phi,
   \: \phi : C( \sigma(a) ) \longrightarrow {\bf B}$. The isomorphism $\phi$
   is isometric and preserves order relations between self-adjoint elements
   and, hence, positivity of self-adjoint elements. Therefore, $\phi$ is a
   positive mapping.

   Choosing an arbitrary maximal abelian $C^*$-subalgebra {\bf D} of {\bf A}
   containing {\bf B} one can complete {\bf B} to $\hat{{\bf B}}({\bf D})$
   with respect to order convergence inside {\bf D}. But $\hat{{\bf B}}
   ({\bf D}) \subseteq {\bf A}$ does not depend on the choice of {\bf D}.
   This can be easily seen if one extends the map $\phi$ to an order
   preserving, isometric mapping $\hat{\phi}$ from the set of all bounded
   complex-valued functions on $\sigma(a)$ into $\hat{{\bf B}}({\bf D})$
   like in \cite{Wr1}, \cite{Wr2}. The characteristic functions of the Borel
   sets of $\sigma(a)$ generate all projections of $\hat{{\bf B}}({\bf D})$
   via $\hat{\phi}$. Moreover, every projection $p$ of $\hat{{\bf B}}
   ({\bf D})$ is the supremum of the net $P = \{ a \in {\bf B}_h^+ : a \le p
   \}$, (where the supremum is taken inside $\hat{{\bf B}}({\bf D})$),
   \cite[Lemma 1.7]{Hm3}. It does not depend on the choice of {\bf D} since
   $P = P^2$ and, hence,  the supremum of $P$ inside every other maximal
   commutative $C^*$-algebra {\bf D}' of {\bf A} containing {\bf B} is also a
   projection $p{\rm '}$. The projections $(1_A -p)$ and $(1_A -p')$ both
   annihilate $P$ inside {\bf A}. The latter implies $p = p'$ because of the
   maximality of {\bf D} and {\bf D}' and of the normality of {\bf A}.

   Now, by \cite[Th.4.1]{Wr1} there exists a unique positive quasi-regular
   $\hat{{\bf B}}$-valued measure $m$ with the property that
   \[
   \int_{\sigma(a)} f(\lambda) \:dm_\lambda \quad = \quad \phi(f)
   \]
   for every $f \in C( \sigma(a)) $. Since $\phi^{-1}(a)(\lambda) = \lambda$
   for every $\lambda \in  \sigma(a) \subset {\bf C}$ by the definition of
   $\phi$ one gets
   \[
   \int_{\sigma(a)} \lambda \: dm_\lambda \quad = \quad a .
   \]
   Moreover, since $\hat{\phi}(\chi_E)^2 =\hat{\phi}(\chi_E^2)=\hat{\phi}
   (\chi_E)  $ for the characteristic function of every Borel set $E \in
   \sigma(a)$ the measure $m$ is projection valued.        $\: \bullet$

   \medskip
   The following corollary is the key point for the subsequent considerations:

   \medskip
   {\bf Corollary 3.3 :} {\it Let} {\bf A} {\it be an} $AW^*${\it -algebra and}
   $\{ \cal M, \langle .,. \rangle \}$ {\it be a pre-Hilbert} {\bf A}{\it
   -module. If} $\: x \in \cal M$ {\it is different from zero then there
   exists a projection} $p \in {\bf A}_h^+$, $p \not= 0$, {\it and an element}
   $a \in {\bf A}_h^+$ {\it such that} $\: a$, $p$ {\it and} $\langle x,x
   \rangle^{1/2}$ {\it commute pairwise, and such that}
   $a \langle x,x \rangle^{1/2} = \langle ax , ax \rangle^{1/2} = p$.

   \medskip
   P r o o f: Consider the commutative $C^*$-subalgebra {\bf B} of {\bf A}
   generated by the elements $\{ 1_A, \langle x,x \rangle \}$. By
   Theorem 3.2 there exists a unique positive quasi-regular measure $m$ on
   the Borel sets of  $\sigma (\langle x,x \rangle^{1/2}) \subset {\bf R}^+$
   being projection-valued in the monotone closure $\hat{{\bf B}}({\bf D})$
   of {\bf B} with respect to an arbitrarily fixed, maximal commutative
   $C^*$-subalgebra {\bf D} of {\bf A} which contains {\bf B}, and satisfying
   the equality
   \[
   \int_{ \sigma( \langle x,x \rangle^{1/2} ) } \lambda \: dm_\lambda \quad =
   \quad \langle x,x \rangle^{1/2}
   \]
   in the sense of order convergence in $\hat{{\bf B}}({\bf D}) \subseteq
   {\bf A }$. Now, if $\langle x,x \rangle^{1/2}$ is a projection set
   \linebreak[4] $a= 1_A$, $p=  \langle x,x \rangle$. If $\langle x,x
   \rangle^{1/2}$ is invertible in {\bf A} set $p= 1_A$, $a= \langle x,x
   \rangle^{-1/2}$. Other\-wise consider a number $\mu \in \sigma(\langle x,x
   \rangle^{1/2})$, $ 0 < \mu < \|x\|$, and set $K=[0,\mu] \cap \sigma
   (\langle x,x \rangle^{1/2})$. The value $m(K) \in \hat{{\bf B}}({\bf D})$
   is a projection different from zero. It commutes with every spectral
   projection of $\langle x,x \rangle^{1/2}$ and with $\langle x,x
   \rangle^{1/2}$ itself. Since $m$ is quasi-regular one has
   \[
   \int_{\sigma(\langle x,x \rangle^{1/2}) \setminus K} \lambda \:
   d(m_\lambda(1_A-m(K)) \: = \: (1_A-m(K)) \langle x,x \rangle^{1/2}.
   \]
   Therefore, one finds $p=(1_A-m(K))$ and $a=((1_A-m(K)) \langle x,x
   \rangle^{-1/2}$, where the inverse is taken in the $C^*$-subalgebra
   $(1_A-m(K)) \hat{{\bf B}}({\bf D}) \subseteq {\bf A}$. Since $\mu < \|x\|$
   the projection $p$ is different from zero. The existence of $a \in
   {\bf A}^+_h $ is guaranteed by $0<\mu$.  $\: \bullet$

   \medskip
   For the completeness of the current section we show that polar
   decomposition is possible inside every $AW^*$-algebra. This generalizes
   assertions of {\sc S.~K.~Berberian} \cite[\S 21, Prop.1,
   Prop.2, Exerc.1,2]{Ber}, {\sc I.~Kaplansky} \cite[Th. 65]{Kapl:68} and
   {\sc R.~V.~Kadison}, {\sc G.~K.~Pedersen}
   \cite[Prop. 2.3]{KP}.

   \medskip
   {\bf Proposition 3.4 :} {\it Let} {\bf A} {\it be an} $AW^*${\it -algebra.
   For every} $x \in {\bf A}$ {\it there exists a unique partial isometry}
   $u \in {\bf A}$ {\it such that} $x= (xx^*)^{1/2} u$ {\it and such that}
   $u u^*$ {\it is the range projection of} $(xx^*)^{1/2}$.

   \medskip
   This follows from Corollary 3.3 and from the above cited results of
   {\sc S.~K.~Berberian}, \linebreak[4] {\sc I.~Kaplansky} 
   and {\sc R.~V.~Kadison}, {\sc G.~K.~Pedersen}.

\section{A criterion on self-duality and $C^*$-reflexivity }

   {\bf Theorem 4.1 :} (cf.~\cite[Th.3.2]{Fr1}) {\it Let} {\bf A} {\it be a
   monotone complete} $C^*${\it -algebra and} $\{ \cal M, \langle .,. \rangle
   \}$ {\it a Hilbert} {\bf A}{\it -module. The following conditions are
   equivalent}:

   (i)  $\;\: \cal M$ {\it is self-dual}.

   (ii) $\: \cal M$ {\it is} $C^*${\it -reflexive}.

   (iii)  {\it The unit ball of } $\cal M$ {\it is complete with respect to}
   $\tau_1^0${\it -convergence}.

   (iv) $\,${\it The unit ball of} $\cal M$ {\it is complete with respect to}
   $\tau_2^0${\it -convergence}. \newline
   {\it If} {\bf A} {\it is commutative there is a further equivalent
   condition, (cf}. \cite[Th.5.4]{Oz1}):

   (v) $\:\; \cal M $ {\it is a Kaplansky-Hilbert module over} {\bf A}.

   \medskip
   Recall the definition of Kaplansky-Hilbert modules
   over commutative $AW^*$-algebras:

   \medskip
   {\bf Definition 4.2:} (\cite[p.842, Def.]{Kp2})
   Let {\bf A} be a commutative $AW^*$-algebra. A Hilbert {\bf A}-module
   $\{ \cal M, \langle .,. \rangle \}$ is {\it Kaplansky-Hilbert} if and
   only if it has the following two properties:

   (i) $\:$ Let $\{ p_\alpha : \alpha \in I \}$ be a set of pairwise
   orthogonal projections of {\bf A} with least upper bound $p \in {\bf A}$.
   Let $x \in \cal M$ be an element for which $p_\alpha x = 0$ for every
   $\alpha \in I$. Then $px=0$.

   (ii) Let $\{ p_\alpha : \alpha \in I \}$ be a net of pairwise orthogonal
   projections of {\bf A} and let
   $\{ x_\alpha : \alpha \in I \}$ be any bounded set in $\cal M$. Then
   there exists an element $x \in \cal M$ such that $p_\alpha x = p_\alpha
   x_\alpha$ for every $\alpha \in I$.

   \medskip
   Before we start proving the theorem we make a simple observation:

   \medskip
   {\bf Lemma 4.3 :} {\it Let} {\bf A } {\it be a monotone complete}
   $C^*${\it -algebra and} $ \{ \cal M, \langle .,. \rangle \}$ {\it be a
   Hilbert} {\bf A}{\it -module possessing a} $\tau_1^0${\it -complete unit
   ball. If} $f: \; \cal M \longrightarrow {\bf A}$ {\it is an} {\bf A}{\it
   -linear bounded mapping for which the set} ${\rm Ker}(f)^\bot = \{ x \in
   \cal M$ : $ \langle x,y \rangle =0$ {\it for every} $y \in {\rm Ker}(f) \}$
   {\it consists only of the zero element then} $f \equiv 0$ {\it on}
   $\cal M$.

   \medskip
   P r o o f: Suppose there exists an element $y \in {\cal M} \setminus
   {\rm Ker}(f)$ with the property $f(y) \neq 0$. One has to show that there
   exist elements $p,a \in {\bf A}^+_h$ for $y \in \cal M$ with the properties
   described in Corollary 3.4 (i.e., $p= \langle ay,ay \rangle$, in
   particular), and with $f(ay) \neq 0$. Indeed, if $f(ay)=0$ for every
   possible choice of $p$ and $a$ then by the equality $f(ay)=(pap)(pf(y))=0$
   and by the invertibility of $a$ inside $p{\bf A}p$ one has $pf(y)=0$.
   Consequently, $f(y)^* p f(y) = 0$ for every possible chosen $p \in
   {\bf A}^+_h$. But the supremum of all such projections $p$ is the support
   of $y \in \cal M$ (cf.~Cor.~3.4), and by \cite[Lemma 1.9]{Hm2} $f(y)=0$
   follows in contradiction to our choice at the beginning. Now denote by
   ${\rm L}_A (ay)$ the {\bf A}-submodule of $\cal M$ generated by the
   element $ay \in \cal M$. Since ${\rm L}_A (ay)$ is generated
   by a single element and norm-closed it is self-dual as a Hilbert
   {\bf A}-submodule of $\cal M$, (\cite[Cor.]{Mi}). Hence, there exists a
   non-zero element $z \in {\rm L}_A (ay)$ such that
   \[
      f|_{{\rm L}_A (ay)} \:\: (x) = \langle x,z \rangle
   \]
   for every $x \in {\rm L}_A (ay)$. But $z \in {\rm Ker}(f)^\bot$ in
   contradiction to the assumption. $\: \bullet$

   \medskip
   P r o o f of the theorem: First, we prove the implication (iii)
   $\rightarrow$ (i). If the Hilbert {\bf A}-module $\cal M$ has a
   $\tau_1^0$-complete unit ball then we can suppose $\langle \cal M, \cal M
   \rangle ={\bf A}$ without loss of generality. Indeed, $\langle \cal M,
   \cal M \rangle$ is a two-sided, monotone and norm-closed $C^*$-ideal of
   {\bf A}. Thus, by \cite[Cor.2.3.1]{Kp2} a central projection $p \in {\bf A},
   p \not= 0$, exists for which $\langle \cal M, \cal M \rangle =$$ {\bf A}p$.
   Therefore, the Hilbert {\bf A}-module $\cal M$ could be considered as a
   Hilbert ${\bf A}p$-module replacing {\bf A} by {\bf A}{\it p}.

   Consider an arbitrary {\bf A}-linear, bounded mapping $f:\cal M
   \longrightarrow {\bf A}$. By Lemma 4.3 one can suppose Ker$(f)^\bot \not=
   \{ 0 \}$. The pair $\{ {\rm Ker}(f)^\bot , \langle .,. \rangle \}$ defines
   a Hilbert {\bf A}-submodule of $\cal M$ possessing a $\tau_1^0$-complete
   unit ball. The map $f$ is faithful on Ker$(f)^\bot$. The image
   $f({\rm Ker}(f)^\bot)$ is a norm-closed left ideal $I$ of {\bf A}, and
   there exists a projection $p \in {\bf A}_h^+$ such that $I={\bf A}p$.
   Indeed, since $\{ {\rm Ker}(f)^\bot , \langle .,. \rangle \}$ and $\{ I,
   \langle .,. \rangle_A \}$ are isomorphic as Hilbert {\bf A}-modules via
   $f$ by assumption one has
   \[
   \langle f(y),f(y) \rangle_A \:\le\: \|f\|^2 \langle y,y \rangle ,  \quad
   \langle f^{-1}(a), f^{-1}(a) \rangle \:\le\: \|f^{-1}\|^2 \langle a,a
   \rangle_A
   \]
   for every $a \in I, y \in {\rm Ker}(f)^\bot$ , (because of the
   {\bf A}-linearity of $f$). But Ker$(f)^\bot$ has a $\tau_1^0$-complete
   unit ball. Therefore, the inequalities above yield $I={\bf A}p$ for a
   projection $p \in {\bf A}^+_h$ as desired since $\|f^{-1}\|$ is
   bounded.

   The Hilbert {\bf A}-module $\{ {\bf A}p, \langle .,. \rangle_A \}$ is
   self-dual, and so is $\{ {\rm Ker}(f)^\bot, \langle .,. \rangle \}$.
   Consequently, there exists an element $x \in {\rm Ker}(f)^\bot$ such that
   $f(.) \equiv \langle .,x \rangle $ on Ker$(f)^\bot$ . Since Ker$(f)$ has
   a $\tau_1^0$-complete unit ball, (i.e., Ker$(f)=({\rm Ker}(f)^\bot)^\bot$),
   the Hilbert {\bf A}-submodule \linebreak[4] $\cal N=$ Ker$(f)+
   $Ker$(f)^\bot$ of $\cal M$ possesses a $\tau_1^0$-complete unit ball. The
   {\bf A}-linear, bounded map $\: g(.)=f(.)-\langle .,x \rangle \:$ is
   obviously equal to zero on $\cal N$ and $\cal N^\bot = $$\{ 0 \}
   \supseteq$ Ker$(g)^\bot$. Therefore, by Lemma 4.3 $g \equiv 0$ on $\cal M$,
   which proves the implication.

   Now we show the implications (i) $\rightarrow$ (ii) $\rightarrow$ (iii).
   The first implication is trivial by definition so we concentrate our
   attention to the second one. Suppose a $C^*$-reflexive Hilbert
   {\bf A}-module $\{ \cal M, \langle .,. \rangle \}$  does not have a
   $\tau_1^0$-complete unit ball. Then consider the Hilbert {\bf A}-module
   $\{ {\cal M}^* , \langle .,. \rangle_D \}$ being derived from $\cal M$ in
   the way described at Proposition 2.6. Since $\cal M^*$ is self-dual by the
   assertions above, since $\cal M \hookrightarrow \cal M^*$ and since
   $\cal M^* \hookrightarrow \cal M'$, the Banach {\bf A}-module $\cal M=
   \cal M''$ has to be identical with $\cal M^*$ in contradiction to
   our assumption.

   The conditions (i) and (iv) are obviously equivalent. Indeed, if a
   non-self-dual Hilbert {\bf A}-module would possess a $\tau_2^0$-complete
   unit ball then its unit ball would have to be non-$\tau_1^0$-complete by
   the previous observations. But this contradicts its $\tau_2^0$-completeness
   by Lemma 2.5. In the same manner one shows the inverse implication.

   In the case of {\bf A} being commutative the equivalence of (v) with some
   of the conditions was shown by {\sc M.~Ozawa} \cite{Oz1}. For completeness,
   we give another proof using methods of order convergence. Let $\{ p_\alpha
   : \alpha \in I \}$ be a set of pairwise orthogonal projections of {\bf A}
   possessing a least upper bound $p \in {\bf A}$. Let $x \in \cal M$ such
   that $p_\alpha x=0$ for every $\alpha \in I$. Because of the inequality
   \[
   0 \:=\: p_\alpha \langle x,x \rangle p_\alpha \le p \langle x,x \rangle p
   \]
   being valid for every $\alpha \in I$, and because of the equality
   $0=\sup \{ p_\alpha \langle x,x \rangle p_\alpha \} = p \langle x,x \rangle
   p$ being valid by (iii) and by the commutativity of {\bf A}  one has $px=0$.
   So one has item (i) of Definition 4.2. Furthermore, let $p=1_A$ and
   let $\{ x_\alpha : \alpha \in I \}$ be a bounded set in $\cal M$ indexed
   by $I$. If $\cal F$ is the net of all finite subsets of $I$ partially
   ordered by inclusion then
   \[
   x\:=\: \tau_1^0-\lim \bigg\{ \sum_{i \in S} p_i x_i \: : \: s \in \cal F
   \bigg\} \quad \in \cal M
   \]  
   and $x$ has the properties to satisfy Definition 4.2,(ii) because of (iii).
   Therefore (iii) implies (v). That (v)  yields (i) was proved by
   {\sc I.~Kaplansky} \cite[Th.5]{Kp2}. So we are done.  $\: \bullet$

   \medskip
   {\bf Examples 4.4 :} (a) Let $X$ be a stonean space and let $C(X)={\bf A}$
   be the $AW^*$-algebra of all continuous, complex valued functions on $X$.
   Suppose $X$ consists of infinitely many points and, therefore, contains at
   least one accumulation point $x \in X$. Let $C_o(X)$ be the set of all
   functions of $C(X)$ vanishing at $x$. Obviously, $C_o(X)$ is an ideal in
   $C(X)$. Setting $\cal M=$$C_o(X)$ and $\cal N=$$C(X)$ with the
   {\bf A}-valued inner product $\langle .,. \rangle_A$ one has
   $\cal M'=\cal M''=\cal N=$$C(X)$. That is, $\cal M$ is not self-dual.

   (b) Let {\bf A} be a monotone complete $C^*$-algebra and consider the
   standard coun\-tably generated Hilbert {\bf A}-module
   \[
   l_2({\bf A})= \bigg\{ a=\{ a_i : i \in {\bf N} \}: a_i \in {\bf A},
   \sum^\infty_{i=1} a_i a_i^* \quad {\rm is} \:\; {\bf A}{\rm
   -norm-convergent} \bigg\},
   \]
   \[
   \langle a,b \rangle = \sum^\infty_{i=1} a_i b_i^*   .
   \]
   It is self-dual if and only if {\bf A} is finite dimensional as a linear
   space, (\cite[Th.4.3]{Fr1}). The Banach {\bf A}-module
   \[
   l_2({\bf A}')= \bigg\{ a=\{ a_i : i \in {\bf N} \} : a_i \in {\bf A}, \sup
   \big\{ \big\| \sum_{i=1}^N a_i a_i^* \big\| : N \in {\bf N} \big\} <
   \infty \bigg\}
   \]
   turns into a self-dual Hilbert {\bf A}-module if one defines the
   {\bf A}-valued
   inner product by
   \[
   \langle a,b \rangle = {\rm LIM} \bigg\{ \sum_{k=1}^N a_k b_k^* : N \in
   {\bf N} \bigg\} .
   \]

   \medskip
   {\bf Corollary 4.5 :} {\it Let} {\bf A} {\it be a monotone complete}
   $C^*${\it -algebra and} $\{ \cal M, \langle .,. \rangle \}$ {\it be a
   self-dual Hilbert} {\bf A}{\it -module such that} $a=0 \in {\bf A}$ {\it
   is the only element of} {\bf A} {\it for which } $a \cal M = $$\{ 0 \}$.
   {\it Then there exists an element} $z \in \cal M$ {\it with the property}
   $\langle z,z \rangle = 1_A$ {\it and} $\langle \cal M, \cal M \rangle =
   {\bf A}$.

   \medskip
   P r o o f: 
   One can choose a maximal set $\{ x_\alpha : \alpha \in I \}$ of
   elements of $\cal M$ with respect to the conditions: (i) $\langle x_\alpha
   , x_\alpha \rangle =p_\alpha =p_\alpha^2 \not= 0$,
   (ii) $ p_\alpha p_\beta = 0$ for every $ \alpha \not= \beta$ since the unit
   ball of $\cal M$ is $\tau_1^0$-complete, (cf. Corollary 3.3
   and {\sc Zorn}'s lemma). Since the equality \linebreak[4] $\sup \{ p_\alpha
   : \alpha \in I \} = 1_A$
   is valid one can define
   \[
   z \:=\: \tau_1^0-\lim \bigg\{  \sum_{i \in S} x_i : S \in \cal F \bigg\}
   \]
   (where $\cal F$ is the net of all finite subsets of $I$),and one finds the
   desired element $z \in \cal M$, and $\{ \langle az,bz \rangle : a,b \in
   {\bf A} \} = {\bf A} \subseteq  \langle \cal M, \cal M \rangle$.
   $\: \bullet$

   \medskip
   The following corollary generalizes a proposition of {\sc W.~L.~Paschke}
   \cite[Prop.3.11]{Pa1},
   Corollary 3.3 and Proposition 3.4. It shows that one has something like
   polar decomposition inside self-dual Hilbert {\bf A}-modules over monotone
   complete $C^*$-algebras {\bf A}.

   \medskip
   {\bf Corollary 4.6 :}  {\it Let} $\cal M$ {\it be a self-dual Hilbert}
   {\bf A}{\it -module over a monotone complete} $C^*${\it -algebra} {\bf A}.
   {\it Every} $x \in \cal M$ {\it can be decomposed} $x=\langle x,x
   \rangle^{1/2} u$, {\it where} $u \in \cal M$ {\it is such that} $\langle
   u,u \rangle$ {\it is the range projection of} $\langle x,x \rangle^{1/2}$.
   {\it This decomposition is unique in the sense that if the equality }
   $x=bv$ {\it is valid for} $b \in {\bf A}^+_h$, $v \in \cal M$ {\it such
   that} $ \langle v,v \rangle$ {\it is the range projection
   of} $b$, {\it then} $v=u$ {\it and} $b=\langle x,x \rangle^{1/2}$.

   \medskip
   P r o o f:  (cf.~\cite[proof of Prop.3.11]{Pa1}) For a fixed  $x \in
   \cal M$ one sets $x_n = h_n^{-1} x$ with $h_n =(\langle x,x \rangle +1/n
   \cdot 1_A)^{1/2}$, $(n \in {\bf N})$. From the equality $\langle x_n,x_n
   \rangle =(\langle x,x \rangle + 1/n \cdot 1_A)^{-1} \langle x,x \rangle$
   one infers that $\| x_n \| \le 1$ for every $n \in {\bf N}$. Let $y \in
   \cal M$ be a $\tau^0_2$-accumulation point of the sequence $\{ x_n : n \in
   {\bf N} \}$. Since $x= h_n x_n$, $(n \in {\bf N})$,  and  $\| h_n -
   \langle x,x \rangle^{1/2} \| \longrightarrow 0$ for $n \rightarrow
   \infty$ the equality $\langle x,x \rangle^{1/2} y=x$ follows. Denote by
   $p$ the range projection of $\langle x,x \rangle^{1/2}$ in {\bf A}. One has
   \[
   \langle x,x \rangle^{1/2} (p-p \langle y,y \rangle p) \langle x,x
   \rangle^{1/2} =0.
   \]
   Since $\| y \| \le 1$ the element $(p-p \langle y,y \rangle p)$ is positive.
   One obtains the equality
   \[
   \langle x,x \rangle^{1/2} (p-p \langle y,y \rangle p)^{1/2} =0
   \]
   \[
   {\rm i.e.,} \quad p(p-p \langle y,y \rangle p)^{1/2} =0.
   \]
   Hence, $p=p \langle y,y \rangle p$ , and one sets $u=py$. Then
   $\langle x,x \rangle^{1/2} u = x$ and $\langle u,u \rangle = p$ as desired.

   To show the uniqueness of the decomposition suppose $x=bv$ with $b \in
   {\bf A}^+_h$, $v \in \cal M$ such that $\langle v,v \rangle$ is the range
   projection of $b$. Then $\langle x,x \rangle =b^2$, and $b= \langle x,x
   \rangle^{1/2}$. Also $\langle v,v \rangle =p$. The equality $\langle v-pv,
   v-pv \rangle =0$ forces $v=pv$. Also, $\langle  x,u \rangle = \langle x,x
   \rangle^{1/2}$ at one side and $\langle x,u \rangle = \langle x,x
   \rangle^{1/2} \langle v,u \rangle$ at the other side. Thus, $\langle x,x
   \rangle^{1/2} (p- \langle v,u \rangle)=0$, and
   \[
   0= p(p- \langle v,u \rangle) = p- \langle pv,u \rangle = p- \langle v,u
   \rangle.
   \]
   Hence, $\langle u-v,u-v \rangle = 0$ and $u=v$.  $\: \bullet$

   \medskip
    Now we get the following general solution of {\sc W.~L.~Paschke}'s problem
    claimed in the introduction:

    \medskip
    {\bf  Theorem 4.7 :} {\it Let {\bf A} be a $C^*$-algebra. For every
    Hilbert} {\bf A}{\it -module} $ \{ \cal M, \langle .,. \rangle \}$ {\it
    the} {\bf A}{\it -valued inner product} $\langle .,. \rangle$ {\it on}
    $\cal M$ {\it can be continued to an} {\bf A}{\it -valued inner product}
    $\langle .,. \rangle_D$ {\it on the} {\bf A}{\it -dual Banach}
    {\bf A}{\it -module} $\cal M'$ {\it turning} $\{ \cal M',\langle .,.
    \rangle_D \}$ {\it into a self-dual Hilbert} {\bf A}{\it -module if and
    only if {\bf A} is monotone complete (iff {\bf A} is additively complete).
    Moreover, the equalities}
    \[
    \langle x,y \rangle_D = \langle x,y \rangle \:,\: \langle x,r \rangle_D
    = r(x)
    \]
    {\it are satisfied for every} $x,y \in \cal M \hookrightarrow \cal M'$,
    {\it every} $r \in \cal M'$.

    \medskip
    P r o o f : One direction follows immediately either from {\sc M.~Hamana}
    \cite[Th. 2.2]{Ham92} or from {\sc H.~Lin} \cite[Lemma 3.7]{Lin} or from
    Proposition 2.6 and Theorem 4.1 above. The converse can be seen combining
    the result of {\sc M.~Hamana} \cite{Ham92} that the $C^*$-algebra {\bf A}
    has to be additively complete with the result of
    {\sc K.~Sait{\^o}} and {\sc J.~D.~M.~Wright} \cite[\S 3]{Sai/Wri} that
    additively complete C*-algebras are monotone complete, and vice versa.
    $\: \bullet$

    \medskip
    {\bf Corollary 4.8 :} {\it Let} {\bf A} {\it be a monotone complete}
    $C^*${\it -algebra and} $ \cal M$ {\it be  a self-dual Hilbert}
    {\bf A}{\it -module. For every pre-Hilbert} {\bf A}{\it -submodule}
    $\cal N \subseteq \cal M$ {\it one can decompose} $\cal M$ {\it into the
    direct sum of} $\, {\cal N}^* = \cal N'$ {\it and} $\cal (N^*)^\bot =
    (\cal N')^\bot$.

    \medskip
    This is the consequence of Proposition 2.6, Theorem 4.7, Theorem 4.1
    and \cite[Th.2.8]{Fr1}.

\section{A structural criterion on self-duality and $C^*$-reflexivity}

    In the present section we want to show a structural criterion. It was
    suggested by \cite[Th.3.11]{Pa1} and first proved by {\sc M.~Hamana}
    \cite[Th. 1.2]{Ham92} independently. We use our own methods to give
    another proof of it. To formulate the assertion we need the following
    definition:

    \medskip
    {\bf Definition 5.1 :} Let $I$ be an index set and let $\{ \{
    {\cal M}_\alpha , \langle .,. \rangle_\alpha \}: \alpha \in I \} $
    be a set of pre-Hilbert {\bf A}-modules over a fixed monotone complete
    $C^*$-algebra {\bf A}. Let $ \cal F$ be the net of all finite subsets of
    $I$ partially ordered by inclusion. Define the value $\langle x,y
    \rangle_S \in {\bf A}$ for all $I$-tuples $x=\{ x_\alpha \in
    {\cal M}_\alpha : \alpha \in I \}$, $y= \{ y_\alpha \in {\cal M}_\alpha :
    \alpha \in I \}$ and for every $S \in \cal F$ by the formula
    \[
    \langle x,y \rangle_S \:=\: \sum_{i \in S} \langle x_i , y_i \rangle_i
    \: .
    \]
    Let $\cal M$ be the set of all $I$-tuples $x=\{ x_\alpha \in
    {\cal M}_\alpha : \alpha \in I \}$ for which the least upper bound
    $\sup \{ \langle x,x \rangle_S : S \in \cal F \}$ exists in {\bf A}. Then
    define for $x,y \in \cal M$
    \[
    \langle x,y \rangle \:=\: {\rm LIM} \{ \langle x,y \rangle_S : S \in
    \cal F \}.
    \]
    The linear space $\cal M$ is a (left) {\bf A}-module with respect to the
    coordinatewise operations induced from the $\{ {\cal M}_\alpha : \alpha
    \in I \}$. Moreover, the mapping $\langle .,. \rangle : {\cal M} \times
    \cal M \longrightarrow {\bf A}$ has all the properties of an
    {\bf A}-valued inner product on $\cal M$ by Lemma 2.2 and Lemma 2.5.
    Let us denote the pre-Hilbert {\bf A}-module $\{ \cal M , \langle .,.
    \rangle \}$ by $\tau^0_1-\Sigma \{ {\cal M}_\alpha : \alpha \in I \}$.
    Note that $\cal M$ is norm-complete (resp., self-dual) if and only if
    each ${\cal M}_\alpha, \alpha \in I$, is.

    \medskip
    {\bf Theorem 5.2 :} {\it Let} {\bf A} {\it be a monotone complete}
    $C^*${\it -algebra and let} $\{ {\cal M}, \langle .,. \rangle \}$ {\it
    be a Hilbert} {\bf A}{\it -module}. {\it The following two conditions are
    equivalent}:

    (i) $\; \cal M$ {\it is self-dual}.

    (ii) {\it There exists a set} $\{ p_\alpha : \alpha \in I \}$ {\it of not
    necessarily distinct projections of} {\bf A} {\it such that} $\cal M$
    {\it and} $\tau_1^0-\Sigma \{ {\bf A} p_\alpha : \alpha \in I \}$ {\it
    are isomorphic as Hilbert} {\bf A}{\it -modules}.

    \medskip
    P r o o f: The implication (ii)$\rightarrow$(i) easily follows from
    Theorem 4.1 and from the self-duality of the Hilbert {\bf A}-modules
    $\{ \{ {\bf A} p_\alpha , \langle .,. \rangle_A \} : \alpha \in I \}$.
    There remains to show the converse.

    By Corollary 3.4 and {\sc Zorn}'s lemma one finds a maximal set of elements
    of $\cal M$ , \linebreak[4] $\{ x_\alpha : \alpha \in I \}$, with respect
    to the assumptions
    (a) $\langle x_\alpha, x_\alpha \rangle = p_\alpha = p^2_\alpha \not= 0$,
    (b) $\langle x_\alpha, x_\beta \rangle = 0$ if $\alpha \not= \beta$.
    Let $\cal F$ be the net of all finite subsets of $I$ being partially
    ordered by inclusion. The equality $x_\alpha = p_\alpha x_\alpha$ is
    valid for every $\alpha \in I$ and, thus, one can define a mapping
    \[
    T : {\cal M} \longrightarrow \tau^0_1-\Sigma \{ {\bf A} p_\alpha : \alpha
    \in I \} \:, \:
    T(x) = \{ \langle x, x_\alpha \rangle : \alpha \in I \}.
    \]
    This mapping $T$ is obviously {\bf A}-linear. To show the surjectivity of
    $T$ consider  \linebreak[4] $\{ a_\alpha p_\alpha : \alpha \in I \} \in
    \tau^0_1-\Sigma \{ {\bf A} p_\alpha : \alpha \in I \}$. Define $y_S = \{
    a_i x_i : i \in S \} \in \cal M$ for every $S \in \cal F$ and $y=
    \tau_2^0-\lim \{ y_S : S \in \cal F \}$, (cf. Th.4.1). One has $y \in
    \cal M$ and $T(y) = \{ a_\alpha p_\alpha : \alpha \in I \}$.
    To prove that $T$ is one-to-one suppose the existence of a non-zero
    element $x \in \cal M$ for which $ \langle x,x_\alpha \rangle = 0$ holds
    for every $\alpha \in I$. By Corollary 3.4 the equality $\langle ax,ax
    \rangle =p=p^2 \not= 0$ holds for some $p \in {\bf A}_h, a \in
    {\bf A}_h^+$. Moreover, since $\langle ax,x_\alpha \rangle = 0$ for every
    $\alpha \in I$ and since the set $\{ x_\alpha : \alpha \in I \}$ is chosen
    to be maximal with respect to (a), (b) one gets $ax=0$ and $p=0$ in
    contradiction to the choice of $p$.
    Finally, one has to show that $\langle T(x),T(x) \rangle = \langle x,x
    \rangle$ holds for every $x \in \cal M$. Indeed, for every $x \in \cal M$
    and every $S \in \cal F$,
    \[
    \langle x_S, x_S \rangle = \sum_{i \in S} \langle x,x_i \rangle p_i
    \langle x_i,x \rangle = \langle T(x_S),T(x_S) \rangle
    \]
    where $x_S = \Sigma \{ \langle x,x_i \rangle x_i : i \in S \}$ by
    definition. The desired equality now obtains by taking the
    $\tau^0_1$-limit on both sides of this equality.  $\: \bullet$

\section{Applications}

    First, we formulate a classification of self-dual (and hence,
    $C^*$-reflexive), countably gene\-rated Hilbert $C^*$-modules over
    monotone complete $C^*$-algebras.

    \medskip
    {\bf Theorem 6.1 :} {\it Let} {\bf A} {\it be a monotone complete}
    $C^*${\it -algebra and let} $\cal M$ {\it be a self-dual, countably
    generated Hilbert} {\bf A}{\it -module. Then there exist only the
    following two possibilities for the inner structure of} {\bf A} {\it and
    of} $\cal M$:

    (i) $\: \cal M$ {\it is finitely generated and } {\bf A} {\it is
    arbitrary}.

    (ii) $\cal M$ {\it is decomposable into the direct sum of a finitely
    generated Hilbert} {\bf A}{\it -module and a countably generated Hilbert}
    {\bf B}{\it -module, where} {\bf B} {\it is a finite-dimensional,
    two-sided} $C^*${\it -ideal of} {\bf A}.

    \medskip
    P r o o f: The first statement follows from \cite[Cor.]{Mi}. One has to
    show the second one and the completeness of the classification. Suppose,
    the Hilbert {\bf A}-module $\cal M$ is self-dual and countably generated.
    By Theorem 5.2 one has
    \[
    {\cal M} = \tau^0_1-\Sigma \{ {\bf A}p_i : i \in {\bf N} \}
    \]
    for some countable set $\{ p_i : i \in {\bf N} \}$ of projections
    of {\bf A}. By an inductive process (dividing in direct summands and
    taking direct sums) one can reach a situation in which the pairwise
    product $p_i p_j = r$ of every two projections $p_i, p_j$ $(i < j)$ of our
    choice is a projection if and only if $r=p_j \not= 0$. Suppose this
    situation is realized. Since $\cal M$ is countably generated the sequence
    \[
    \bigg\{ \sum_{i=1}^N a_i a_i^* : \, N \in {\bf N} \bigg\}
    \]
    has to converge with respect to the {\bf A}-norm for every element
    $a= \{ a_i \in {\bf A}p_i ,    \, i \in {\bf N} \}$ of
    $\tau^0_1-\Sigma \{ {\bf A}p_i : \, i \in {\bf N} \}$. Therefore, if there
    are more than a finite number of infinite dimensional, two-sided
    $AW^*$-ideals $\langle {\bf A}p_i , {\bf A}p_i \rangle$ of our choice,
    then the Hilbert {\bf A}-module $\tau^0_1-\Sigma \{ {\bf A}p_i : \, i \in
    {\bf N} \}$ can not be countably generated, cf.~\cite[Th.4.3]{Fr1}.
    Moreover, if there does not exist a finite dimensional, two-sided
    $C^*$-ideal {\bf B} in {\bf A} containing all the finite dimensional,
    two-sided $C^*$-ideals $\langle {\bf A}p_i , {\bf A}p_i \rangle$ of our
    choice then the Hilbert {\bf A}-module $\tau^0_1-\Sigma \{ {\bf A}p_i :
    \, i \in {\bf N} \}$ can not be countably generated by \cite[Th.4.3]{Fr1},
    again. So the statements follow.  $\: \bullet$

    \medskip
    Secondly, extending \cite[Th.3.7]{Pa1} we show how {\bf A}-linear,
    bounded operators on a Hilbert
    {\bf A}-module $\cal M$ over a monotone complete $C^*$-algebra {\bf A}
    can be continued to {\bf A}-linear, bounded operators on the {\bf A}-dual
    Hilbert {\bf A}-module $\cal M'$ in a unique way.

    \medskip
    {\bf Proposition 6.2:} {\it Let} {\bf A} {\it be a monotone complete}
    $C^*${\it -algebra and let} $\{ {\cal M}, \langle .,. \rangle \}$ {\it be
    a Hilbert} {\bf A}{\it -module. Then every} {\bf A}{\it -linear bounded
    operator} $T : \cal M \longrightarrow \cal M$ {\it can be continued to a
    unique} {\bf A}{\it -linear bounded operator} $T' : \cal M'
    \longrightarrow \cal M'$ {\it on the} {\bf A}{\it -dual Banach}
    {\bf A}{\it -module} $\cal M'$ {\it of} $\cal M$ {\it preserving the
    operator norm. Moreover, if the operator} $T$ {\it has an adjoint
    operator} $T^* : \cal M \longrightarrow \cal M$ {\it then} $(T^*)' = (
    T')^*$.

    \medskip
    P r o o f: By \cite[Th. 2.8]{Pa1} one obtains
    \begin{equation}
    \langle T(x) , T(x) \rangle \leq \|T\| \langle x,x \rangle
    \end{equation}
    for every $x \in \cal M$. Hence, one can define the operator $T' :
    \cal M' \longrightarrow \cal M'$ by the formulae:
    \[
    T'(x) = T(x)  \,\: {\rm for} \: {\rm every} \,\: x \in {\cal M},
    \]
    \[
    T'(\tau^0_1-\lim \{ x_\alpha :\alpha \in I \}) = \tau^0_1-\lim \{
    T(x_\alpha) : \alpha \in I \}
    \]
    where $\{ x_\alpha : \alpha \in I \}$ is an arbitrarily chosen
    $\tau^0_1$-fundamental net of $\cal M$.
    Obviously, the operator $T'$ is {\bf A}-linear and bounded by $\|T\|$,
    cf. Lemma 2.2(iv) and (1). It is unique by (1). And by Lemma 2.5
    one has $(T^*)' = (T')^*$.
    $\: \bullet$

    The following corollary generalizes results of {\sc G.~Wittstock}
    \cite{Wt} :

    \medskip
    {\bf Corollary 6.3 :} {\it Let} {\bf A} {\it be a monotone complete}
    $C^*${\it -algebra, let} $\{ \cal M, \langle .,. \rangle \}$ {\it be a
    self-dual Hilbert} {\bf A}{\it -module and} $\cal N$ {\it be a Hilbert}
    {\bf A}{\it -submodule of} $\cal M$. {\it Then every} {\bf A}{\it -linear
    bounded operator} $T : \cal N \longrightarrow \cal M$ {\it can be
    continued to a unique} {\bf A}{\it -linear bounded operator} $T' : \cal M
    \longrightarrow \cal M$ {\it preserving the operator norm and the relation}
    $T'({\cal N}^\bot) = \{ 0 \}$. {\it The self-duality of $\cal M$ is
    necessary, in general, for the result. }

    \medskip
    P r o o f: Consider the unique extension $T' : \cal N' \longrightarrow
    \cal M'$ with $\|T'\| = \|T\|$ existing by the previous proposition.
    Since $\cal M'=\cal N' \oplus {\cal N}^{\bot}$ one defines the final
    operator on $\cal M$ by $T'$ on $\cal N'$ and by the zero operator on
    ${\cal N}^\bot$.

    To show that the restriction on $\cal M$ to be self-dual can not be
    dropped, in general, one constructs a counterexample. Consider a
    norm-closed left ideal {\bf D} of {\bf A}, where {\bf D} is order dense
    in {\bf A} and unequal to {\bf A}. (For example, take the set of all
    bounded linear operators on a separable Hilbert space as {\bf A} and the
    compact one's as {\bf D}.) Set ${\cal M} = {\bf A} \oplus {\bf D}$ and
    ${\cal N} = {\bf D} \oplus {\bf D}$ with the usual inner products on them.
    Note, that $\cal N$ is contained in $\cal M$ as a submodule. But the
    operator $T: (d_1,d_2) \longrightarrow (0,d_1)$ can not be continued to
    an operator $T'$ in any way. $\: \bullet$

    \medskip
    Thirdly, let $\cal M$ be a self-adjoint Hilbert {\bf A}-module over a
    monotone complete $C^*$-algebra {\bf A}. Denote by ${\bf End}_A(\cal M)$
    the set of all bounded {\bf A}-linear operators on $\cal M$. We consider
    $C^*$-subalgebras {\bf M} of ${\bf End}_A(\cal M)$ coinciding with their
    bicommutant {\bf M}'' inside ${\bf End}_A(\cal M)$. The following fact was
    obtained by {\sc M.~Hamana} \cite[Prop. 1.2]{Ham92} independently: 

    \medskip
    {\bf Theorem 6.4 :} {\it Let} {\bf A} {\it be a monotone complete}
    $C^*${\it -algebra and let} $\cal M$ {\it be a self-dual Hilbert}
    {\bf A}{\it -module. Then every} $C^*${\it -subalgebra} {\bf M} {\it of}
    $\, {\bf End}_A(\cal M)$ {\it which coincides with its bicommutant}
    {\bf M}'' {\it inside} ${\bf End}_A(\cal M)$ {\it is a monotone complete}
    $C^*${\it -algebra}.

    \medskip
    P r o o f: Consider an arbitrary bounded, increasingly directed net
    $\{ B_\alpha : \alpha \in I \}$ of po\-sitive elements of {\bf M}. By
    \cite[Th. 3]{Loy} an element $C \in {\bf End}_A(\cal M)$ is positive if
    and only if \linebreak[4] $\langle C(x),x \rangle \geq 0$ for every $x
    \in \cal M$. Therefore, by Theorem 4.1 there exists an element $B \in
    {\bf End}_A(\cal M)$ being defined by the formula
    \[
    B(x) = \tau^0_1-\lim \{ B_\alpha (x) : \alpha \in I \} , x \in \cal M,
    \]
    and $B = \sup \{ B_\alpha : \alpha \in I \}$ in ${\bf End}_A(\cal M)$.
    But for every $B_\alpha, \alpha \in I, $ and every element $C$ of the
    commutant {\bf M}' of {\bf M} the equality $CB_\alpha - B_\alpha C = 0$
    is valid. Since $BC-CB = {\rm LIM} \{ (CB_\alpha - B_\alpha C) :\alpha
    \in I \} = 0$ inside ${\bf End}_A(\cal M)$ one has $B \in {\bf M}'' =
    {\bf M}$. $\: \bullet$

    \medskip
    This result seems to be of some importance.
    Immediatelly one realizes that polar decomposition and spectral
    decomposition work inside ${\bf M}''={\bf M} \subseteq {\bf End}_A(\cal M)$
    in every such case, (cf.~Theorem 3.2 and 3.6).
    For investigations about similar monotone complete $C^*$-algebras Theorem
    6.4 allows, for example, to introduce a notion of "Morita equivalence in
    order" along the line of the ideas of {\sc M.~A.~Rieffel} \cite{Rie} for
    $W^*$-algebras. Another area of application is the theory of operator
    valued weights and conditional expectations of finite index between
    monotone complete $C^*$-algebras following {\sc M.~Baillet},
    {\sc Y.~Denizeau} and
    \linebreak[4]
    {\sc J.-F.~Havet} \cite{Bai}. But these investigations
    appear elsewhere, cf.~\cite{Fr4}. What we will do is to prove a
    generalized {\sc Weyl}-{\sc Berg} theorem, which will finish up the
    present paper.

 \section{A Weyl-Berg type theorem}

The problem of approximation of normal elements in $C^*$-algebras
{\bf A} by diagonalizable ele\-ments of {\bf A} up to a suitable small
remainder is part of the sphere of interests of many mathe\-maticians. For an
overview on the recent results and open problems compare the papers of
{\sc D.~Voiculescu} \cite{Voi}, {\sc R.~V.~Kadison} \cite{Kd1,Kd2,Kd3},
{\sc K.~Grove} and {\sc G.~K.~Pedersen} \cite{Gro}, {\sc G.~J.~Murphy}
\cite{Mu}, {\sc N.~Higson} and {\sc M.~R{\o}rdam} \cite{HR}, {\sc L.~G.~Brown}
and {\sc G.~K.~Pedersen} \cite{Bro} and {\sc S.~Zhang} \cite{Zh1,Zh2}. What
we would like is to show a {\sc Weyl}-{\sc Berg} type theorem for monotone
complete $C^*$-algebras and some corollaries of it.

\medskip
{\bf Definition 7.1 :} A monotone complete $C^*$-algebra {\bf A} has the
{\it approximation property} (*) if there exists a chain of pairwise
orthogonal projections $\{ p_\alpha : \alpha \in I \}$ of ${\bf A}_h$ with
least upper bound $1_A$ such that for every $\alpha \in I$ the monotone
complete $C^*$-algebra $p_\alpha {\bf A} p_\alpha$ possesses a faithful state
$f_\alpha$ with the property that the norm-completion  of the pre-Hilbert
space $\{ p_\alpha {\bf A} p_\alpha , f_\alpha(\langle .,. \rangle_A) \}$
is separable.

Note, that in the case of {\bf A} being a $W^*$-algebra the states $f_\alpha$
can be chosen as normal states.

\medskip
The class of monotone complete $C^*$-algebras with the approximation property
(*) is sufficiently large to contain most of the $W^*$-algebras of physical
interest. The {\sc Dixmier} algebra D([0,1]) has property (*), too. But there
are remarkable examples of $W^*$-algebras which do not have property (*),
cf.~\cite[Remark after Def. 2.5.1.]{Bra}. Other examples of commutative
$AW^*$-algebras without property (*) were constructed by {\sc M.~Ozawa}
\cite{Oz2} which can be seen by comparing Proposition 7.6 and Theorem 7.3
below.

Let us denote the norm-closure of the linear hull of
\[
\{ \theta_{x,y} \in {\bf End}_A({\cal M}) : \theta_{x,y} (z)= \langle z,x
\rangle y \: {\rm for} \: {\rm every} \: z \in {\cal M}, \: {\rm each} \:
x,y \in \cal M \}
\]
by ${\bf K}_A(\cal M)$ as usually. ${\bf K}_A(\cal M)$ is called the set of
"compact" operators on $\cal M$.

\medskip
{\bf Definition 7.2 :} Let {\bf A} be a monotone complete $C^*$-algebra, let
$\cal M$ be a self-dual Hilbert {\bf A}-module possessing a countably
generated {\bf A}-pre-dual Hilbert {\bf A}-module. Let $0 \not= q=q^2 \in
{\bf A}_h^+ $ be fixed.
An operator $D \in q {\bf End}_A(\cal M)$ is said to be
{\it diagonalizable} if there is a sequence $\{ x_n : n \in {\bf N} \}$
of pairwise orthogonal nontrivial elements of $q {\cal M}$ such that
$D(x_n) = a_n x_n$ for some $\{ a_n \} \in q {\bf A} q$,
the subsets $q{\bf A}x_n=\{ ax_n : a \in q {\bf A} q \}$, $(n \in
{\bf N})$, are norm-closed and $\tau-\Sigma\{ q {\bf A} x_n : n \in {\bf N}
\} = q {\cal M}$.

\medskip
Note that in the case of {\bf A} being commutative {\bf A} can be identified
with the centre of ${\bf End}_A(\cal M)$ and, hence, $D$ is diagonalizable
inside $q{\bf End}_A(\cal M)$. So we are in the classical situation. 

\medskip
{\bf Theorem 7.3 :} {\it Let} {\bf A} {\it be a monotone complete} $C^*${\it
-algebra with property} (*). {\it Let} $\cal M$ {\it be a self-dual Hilbert}
{\bf A}{\it -module possessing a countably generated} {\bf A}{\it -pre-dual
Hilbert} {\bf A}{\it -module. Suppose that $q{\bf End}_A({\cal M}) \not=
q{\bf K}_A({\cal M})$ for every central projection $q \in {\bf A}$.
Then for every} $\varepsilon > 0$ {\it and every self-adjoint operator} $T
\in {\bf End}_A(\cal M)$ {\it there exist a diagonalizable, self-adjoint
operator} $D \in {\bf End}_A(\cal M)$ {\it and a "compact", self-adjoint
operator} $K \in {\bf K}_A(\cal M)$ {\it such that} $T = D + K$ {\it and}
$\| K \| < \varepsilon$. \newline
{\it If $q{\bf End}_A({\cal M}) = q{\bf K}_A({\cal M})$ for a
central projection $q \in {\bf A}$ and if {\bf A} is a $W^*$-algebra or a
commutative $AW^*$-algebra then the Hilbert {\bf A}-module $q \cal M$ is
finitely generated and every self-adjoint bounded module operator is
diagonalizable.
}

\medskip
P r o o f: The $C^*$-algebra {\bf A} has property (*), and \cite[Th.9]{Mu}
is valid. Consequently, for every $\varepsilon > 0$, every $\alpha \in I$
there exist a diagonalizable, self-adjoint operator $D_\alpha \in p_\alpha
{\bf End}_A(\cal M)$ with eigenvalues $a_n \in \{ \lambda 1_A : \lambda \in
{\bf C} \}$ and a "compact", self-adjoint operator $K_\alpha \in p_\alpha
{\bf K}_A(\cal M)$ with $\| K_\alpha \| < \varepsilon$ such that $p_\alpha T
= D_\alpha + K_\alpha$.

Since these operators $\{ D_\alpha \}, \{ K_\alpha \}$ are pairwise orthogonal,
linear operators on $\cal M$ one can sum them up in the sense of order
convergence inside the monotone complete $C^*$-algebra $({\bf A} +
{\bf End}_A(\cal M)) \subset {\bf End}_C(\cal M)$. One gets
\[
T = \sum_{\alpha \in I} p_\alpha T = \sum_{\alpha \in I} D_\alpha +
\sum_{\alpha \in I} K_\alpha = D+K
\]
where $D \in {\bf End}_A(\cal M)$ is diagonalizable and self-adjoint, and
$K \in {\bf K}_A(\cal M)$ is "compact" and self-adjoint by construction.
Moreover, $\| K \| < \varepsilon$.

If $q{\bf End}_A({\cal M}) = q{\bf K}_A({\cal M})$ for a
central projection $q \in {\bf A}$ then the Hilbert {\bf A}-module $q \cal M$
is algebraicly finitely generated by \cite[Remark 4.5]{Skand},
\cite[Prop. 3.2]{Exel}. Since every normal element of ${\bf M}_n({\bf A})$
is diagonalizable for every $n \in \bf N$ and for {\bf A} being a
$W^*$-algebra or a commutative $AW^*$-algebra by \cite{Kd2,Chr/Ped} and since
$q{\cal M})=P(q{\bf A}^n)$ for a natural number $n$ and a
{\bf A}-linear projection $P$ on ${\bf A}^n$ the desired result yields.
$\: \bullet$

\medskip
To extend the statement of Theorem 7.3 to normal elements $T \in
{\bf End}_A(\cal M)$ note that there always exists a self-adjoint element $
S \in {\bf End}_A(\cal M)$ such that $T$ is contained in the $C^*$-subalgebra
being generated by $S$ and the identity of ${\bf End}_A(\cal M)$. This follows
from Theorem 3.1 by functional calculus, (cf. \cite{Hl1,Hl2,Va}).
Now the techniques of {\sc G.~J.~Murphy} \cite[p.283]{Mu} allow to prove the
following:

\medskip
{\bf Theorem 7.4 :} {\it Let} {\bf A} {\it be a monotone complete} $C^*${\it
-algebra with property} (*). {\it Let} $\cal M$ {\it be a self-dual Hilbert}
{\bf A}{\it -module possessing a countably generated} {\bf A}{\it -pre-dual
Hilbert} {\bf A}{\it -module. Suppose that $q{\bf End}_A({\cal M}) \not=
q{\bf K}_A({\cal M})$ for every central projection $q \in {\bf A}$.
Then for every normal operator} $T \in {\bf End}_A(\cal M)$ {\it there exist a
diagonalizable operator} $D \in {\bf End}_A(\cal M)$ {\it and a "compact"
operator} $K \in {\bf End}_A(\cal M)$ {\it such that} $T = D + K$. \newline
{\it If $q{\bf End}_A({\cal M}) = q{\bf K}_A({\cal M})$ for a
central projection $q \in {\bf A}$ and if {\bf A} is a $W^*$-algebra or a
commutative $AW^*$-algebra then the Hilbert {\bf A}-module $q \cal M$ is
finitely generated and every normal bounded module operator is diagonalizable.
}

\medskip
{\bf Remark 7.5 :} The property of $\cal M$ to possess a countably generated
Hilbert {\bf A}-module as it's {\bf A}-pre-dual can not be dropped,
cf.~\cite{Hl2}. But the infinite cardinality of minimal generator sets of some
{\bf A}-pre-dual Hilbert {\bf A}-modules of $\cal M$ has to be unique in some
sense as shown below in Proposition 7.6. Unfortunally, it is not quite clear
at present if {\bf A} has to possess property (*) to validate the
statement of Theorem 7.4, or whether other situations are possible. Beside
this it is strange that a property of {\bf A} should determine the truth of
the generalized {\sc Weyl}-{\sc Berg} theorem with respect to the pair
$\{ {\bf End}_A(\cal M), {\bf K}_A(\cal M) \}$, because the intersection of
the three $C^*$-algebras inside ${\bf End}_C(\cal M)$ is only {\bf Z}({\bf A}),
the centre of {\bf A}. So the translation of (*) to a $C^*$-condition on
$({\bf End}_A(\cal M), {\bf K}_A(\cal M))$ appears to be non-easy if one tries
to formulate it without reference to {\bf A}, (at least if {\bf A} is
non-commutative). However, there are some facts that shead light on the
situation, cf.~Proposition 7.8.

\medskip
To prove the following fact one makes use of a discovery of
{\sc M.~Ozawa} \cite{Oz2}. He showed that the following misbehaviour may happen: 
There are special commutative $AW^*$-algebras {\bf A} such that the self-dual 
Hilbert {\bf A}-module $l_2({\bf A})'$ is isomorphic to the self-dual Hilbert 
{\bf A}-module $\tau^0_1-\Sigma \{ {\bf A}_{(\alpha)} : \alpha \in I \}$ for 
some set $I$ of uncountable cardinality $card(I)$, where $card(I)$ depends
on the inner structure of {\bf A}.

\medskip
{\bf Proposition 7.6 :} {\it Let} {\bf A} {\it be a commutative} $AW^*${\it
-algebra. If for a given real number} \linebreak[4] $\varepsilon > 0$ {\it one
can decompose every self-adjoint operator} $T \in {\bf End}_A(l_2({\bf A})')$
{\it into the sum} $T=D+K$ {\it of a diagonalizable, self-adjoint operator}
$D \in {\bf End}_A(l_2({\bf A})')$ {\it and a "compact", self-adjoint operator}
$K \in {\bf K}_A(l_2({\bf A})')$ {\it with} $\| K \| < \varepsilon$ {\it then
the existence of an {\bf A}-linear isomorphism of} $\: l_2({\bf A})'$ {\it to
a self-dual Hilbert} {\bf A}-module {\it of type} $\tau^0_1-\Sigma \{
{\bf A}_{(\alpha)} : \alpha \in I \}$ {\it for sets} $I$ {\it of uncountable
cardinality} $card(I)$ {\it is impossible, i.e.~the structure of {\bf A} can
not be arbitrary}.

\medskip
P r o o f: Suppose $l_2({\bf A})'$ is isomorphic to $\tau^0_1-\Sigma \{
{\bf A}_{(\alpha)} : \alpha \in I \}$ for some set $I$ of uncountable
cardinality $card(I)$. Denote by $\cal H$ the non-separable Hilbert space
$\tau^0_1-\Sigma \{ {\bf C}_{(\alpha)} : \alpha \in I \}$.  By
\cite{Hl2} there exists a self-adjoint operator $T_0 \in {\bf End}_C(\cal H)$
being not decomposable into the sum of a diagonalizable, self-adjoint operator
and a compact, self-adjoint operator inside ${\bf End}_C(\cal H)$. Since
\[
l_2({\bf A})' \: = \: \tau^0_1-\Sigma \{ {\bf A}_{(\alpha)} : \alpha \in I \}
\: = \: ({\bf A} \otimes \cal H)'
\]
by assumption one can define a self-adjoint, bounded, {\bf A}-linear operator
\newline
$T: ({\bf A} \otimes {\cal H})' \longrightarrow ({\bf A} \otimes \cal H)'$ as
the unique extension of the mapping
\[
T' : \: a \otimes h \longrightarrow a \otimes T_0 (h) \: , \:(a \in {\bf A},
h \in \cal H),
\]
from ${\bf A} \otimes \cal H$ to $({\bf A} \otimes \cal H)'$,
(cf.~Prop.~6.2).
Assume now that $T$ would be decomposable into the sum of a diagonalizable,
self-adjoint operator $D \in {\bf End}_A(l_2({\bf A})')$ and a "compact",
self-adjoint operator $K \in {\bf K}_A(l_2({\bf A})')$. If one considers
{\bf A} as $C(X)$ for some compact Hausdorff space X one finds
\[
T_{(x)} = T_0 = D_{(x)} + K_{(x)}
\]
with $D_{(x)} \in {\bf End}_C({\cal H}_x)$ -- diagonalizable and self-adjoint,
$K_{(x)} \in {\bf K}_C({\cal H}_x)$ -- compact and self-adjoint. But this
contradicts the choice of $T_0$.  $\: \bullet$

\medskip
{\bf Corollary 7.7 :} {\it There are commutative} $AW^*${\it -algebras} {\bf A}
{\it such that every normal \linebreak[4] element of} ${\bf M}_n ({\bf A})$,
($n \in {\bf N}$ {\it -- arbitrary), is diagonalizable, but the generalized
{\sc Weyl}-{\sc Berg} theorem is not valid for some} {\bf A}{\it -linear,
bounded, self-adjoint operators on} $l_2({\bf A})'$.

\medskip
For an example compare \cite{Oz2}, \cite[Added in proof]{Gro} and Proposition
7.5.

\medskip
{\bf Proposition 7.8 :}
Let {\bf A} be a monotone complete $C^*$-algebra possessing a countable chain 
$\{ q_k : k \in {\bf N} \}$ of pairwise orthogonal projections such that
$\Sigma_{k=1}^{\infty} q_k = 1_A$ in the sense of order convergence, and
$q_k \sim 1_A$ for every $k \in {\bf N}$. Let $\cal M$ be a self-dual Hilbert
{\bf A}-module with a countably generated {\bf A}-pre-dual Hilbert
{\bf A}-module. Then ${\bf End}_A({\cal M}) = {\bf K}_A({\cal M})$.

\medskip
P r o o f:  By Theorem 5.2 we have ${\cal M} = \tau^0_1 - \Sigma \{ {\bf A}
p_n : n \in {\bf N} \}$ for projections $\{ p_n \} \in {\bf A}$.
Denote by $u_k$ the partial isometries of {\bf A} realizing the equivalences
$q_k \sim 1_A$, i.e., $u_k u_k^* = q_k$, $u_k^*u_k = 1_A$. Then
\[
id_{\cal M}( \cdot ) = \langle \cdot , u \rangle u = \theta_{u,u} ( \cdot )
\]
where $u = \{ u_k p_k : k \in {\bf N} \} \in \cal M$. That is, the identity
operator on $\cal M$ is ''compact''. Since ${\bf K}_A(\cal M)$ is a
two-sided ideal in ${\bf End}_A(\cal M)$ the proof is complete. $\, \bullet$

\medskip
Similarly, if $\cal M$ is finitely generated then ${\bf End}_A({\cal M}) =
{\bf K}_A({\cal M})$, too, for every unital $C^*$-algebra {\bf A}. Especially,
if {\bf A} is a $W^*$-algebra or a commutative $AW^*$-algebra then every
element of ${\bf End}_A({\cal M})$ is diagonalizable, cf.~\cite{Kd1,Kd2} and
\cite[Cor. 3.3]{DEPE}. Hence, the following two problems are unsolved:

\medskip
{\bf Problem :} Does there exist a countable chain $\{q_k : k \in {\bf N} \}$
of pairwise orthogonal projections inside every monotone complete
$AW^*$-factor of type ${\rm II}_\infty$ or ${\rm III}$ such
that $\Sigma q_k = 1_A$ in the sense of order convergence, and $q_k \sim 1_A$
for every $k \in {\bf N}$ ?

\medskip
{\bf Problem :} Is every normal bounded module operator on finitely generated
Hilbert {\bf A}-modules over monotone complete $C^*$-algebras {\bf A}
diagonalizable?

    \bigskip \noindent
    Universit\"at Leipzig  \newline
    FB Mathematik/Informatik  \newline
    Mathematisches Institut   \newline
    04109 Leipzig          \newline
    Augustusplatz 10            \newline
    Federal Republic of Germany  \newline
    frank@mathematik.uni-leipzig.de

\end{document}